\documentclass[12pt]{amsart}
\font\emailfont=cmtt10

\usepackage{amsmath,amsthm,amsfonts,amscd,flafter,epsf}

\newcommand\commentable[1]{#1}

\newcommand{\rk}{\mathrm{rk}}
\newcommand{\HF}{HF}

\newtheorem{theorem}{Theorem}[section]
\newtheorem{prop}[theorem]{Proposition}
\newtheorem{cor}[theorem]{Corollary}

\newtheorem{lemma}[theorem]{Lemma}

\newtheorem{defn}[theorem]{Definition}

\def\endproof{\relax\ifmmode\expandafter\endproofmath\else
  \unskip\nobreak\hfil\penalty50\hskip.75em\hbox{}\nobreak\hfil\bull
  {\parfillskip=0pt \finalhyphendemerits=0 \bigbreak}\fi}
\def\endproofmath$${\eqno\bull$$\bigbreak}
\def\bull{\vbox{\hrule\hbox{\vrule\kern3pt\vbox{\kern6pt}\kern3pt\vrule}\hrule}}

\newcommand{\Q}{\mathbb{Q}}

\newcommand{\Z}{\mathbb{Z}}

\newcommand{\OneHalf}{\frac{1}{2}}

\newcommand{\OneQuarter}{\frac{1}{4}}

\newcommand{\Zmod}[1]{\Z/{#1}\Z}

\newcommand{\Ker}{\mathrm{Ker}}

\newcommand{\Coker}{\mathrm{Coker}}

\newcommand{\Image}{\mathrm{Im}}

\newcommand{\cm}{\cdot}

\newcommand{\ModSWfour}{\mathcal{M}}
\newcommand{\ModFlow}{\ModSWfour}

\newcommand{\SpinC}{{\mathrm{Spin}}^c}

\newcommand\Hom{\mathrm{Hom}}

\newcommand\abuts\Rightarrow
\newcommand\Sym{\mathrm{Sym}}

\newcommand\mCP{{\overline{\mathbb{CP}}}^2}

\newcommand\HFpRed{\HFp_{\red}}

\newcommand\chiTrunc{\chi^{\mathrm{trunc}}}

\newcommand\relspinc{\underline{\spinc}}

\newcommand\ModSphere{\ModFlow\left({\mathbb S}\longrightarrow 
\Sym^{g-1}(\Sigma_{1})\times \Sym^2(\Sigma_{2})\right)}
\newcommand\ModSpheres\ModSphere

\newcommand{\red}{\mathrm{red}}

\newcommand\HFp{\HFb}

\newcommand\HFinf{HF^\infty}

\newcommand\HFb{HF^+}

\newcommand\UnparModSp{\widehat \ModSp}
\newcommand\UnparModFlow\UnparModSp
\newcommand\Mod\ModSp

\newcommand\PD{\mathrm{PD}}

\newcommand{\spinc}{\mathfrak s}

\newcommand{\spinct}{\mathfrak t}

\newcommand\ModMaps{\mathcal M}
\newcommand\ModSp\ModMaps

\newcommand\Fp[1]{F^{+}_{#1}}

\newcommand\MT{t}

\newcommand\Field{\mathbb F}

\newcommand\HFc{\HF^\circ}

\newcommand\Dual{\mathcal D}
\newcommand\Duality\Dual

\newcommand\FormF[1]{{\mathbb F}^+_{#1}}

\newcommand\FormIFp{{\mathbb H}^+}
\newcommand\MinVecLen{M}

\newcommand\spincrel\relspinc

\commentable{

\title[{Knots with unknotting number one
and Heegaard Floer homology}]
{Knots with unknotting number one
and Heegaard Floer homology}

\author[Peter Ozsv{\'a}th]{Peter Ozsv\'ath}
\address{Department of
Mathematics, Columbia University, New York 10025 \newline
\indent{Institute for Advanced Study, Princeton, New Jersey 08540} \newline
\indent{\emailfont{petero@math.columbia.edu}}}
\thanks{PSO was partially supported by NSF grant numbers DMS-0234311,
DMS-0111298,
and FRG-0244663}

\author[Zolt{\'a}n Szab{\'o}]{Zolt{\'a}n Szab{\'o}} 
\address{Department of
Mathematics, Princeton University, New Jersey 08544 \newline
\indent{\emailfont{szabo@math.princeton.edu}}}}
\thanks{ZSz was partially supported by NSF grant numbers DMS-0107792
and FRG-0244663, and a Packard Fellowship.}
 
\begin{document}

\begin{abstract}  
We use Heegaard Floer homology to give obstructions to unknotting a
knot with a single crossing change.  These restrictions are
particularly useful in the case where the knot in question is
alternating. As an example, we use them to classify all knots with
crossing number less than or equal to nine and unknotting number equal
to one. We also classify alternating knots with ten crossings and
unknotting number equal to one.
\end{abstract}

\maketitle
\section{Introduction}

The {\em unknotting number} $u(K)$ of a knot $K\subset S^3$ is the minimal
number of crossing changes required to unknot  $K$. 
One lower bound on this number is given by a result of
Murasugi~\cite{Murasugi}, 
$${|\sigma(K)|}\leq 2u(K),$$ where here
$\sigma(K)$ denotes the signature of the knot. 

In this article, we focus on obstructions to a knot $K$ having
$u(K)=1$. One classical observation is that if $K$ has $u(K)=1$, then
the branched double-cover of $S^3$ along $K$, $\Sigma(K)$, has cyclic
first homology, c.f. \cite{Nakanishi}. Another obstruction stems from
the linking form of $\Sigma(K)$, c.f.~\cite{Lickorish}.  Indeed, a
number of other obstructions to $K$ having $u(K)=1$ can be given by
considering its branched double-cover.  For example Kanenobu and
Murakami in ~\cite{KanenobuMurakami} use this construction, together
with the cyclic surgery theorem of Culler, Gordon, Luecke, and
Shalen~\cite{CGLS} to classify all two-bridge knots with $u(K)=1$.
More recently, Gordon and Luecke~\cite{GordonLueckeUK1} study
this problem for knots whose branched double-cover is toroidal. 
The aim of the present article is to give a new obstruction for a knot
to have unknotting number equal to one, which uses Heegaard Floer
homology~\cite{HolDisk} for the branched double-cover.

The obstruction is most easily stated in the case where $K$ is an
alternating knot. Indeed, in this case, a statement can be given
purely in terms of elementary number-theoretic properties of the
Goeritz matrix of $K$.  We give some necessary background.

Let 
$$Q\colon V \otimes V \longrightarrow \Z$$ 
be a negative-definite quadratic form over a lattice $V$.
$Q$ determines a map
$$q\colon V \longrightarrow V^*,$$
where here $V^*$ denotes $\Hom(V,\Z)$,
and it induces a bilinear form with values in the rationals
$$Q^*\colon V^* \otimes V^* \longrightarrow \Q.$$
We suppose that the number of elements in the cokernel of $q$ is odd.
A dual element $K\in V^*$ is called {\em
characteristic} if $$\langle K, v\rangle \equiv Q(v\otimes v)\pmod{2}$$ for
all $v\in V$. Fix an element $v_0\in V$ with the property that
$q(v_0)$ is characteristic. We now define 
a map
$$\MinVecLen_Q\colon \Coker(q) \longrightarrow \Q$$
by 
$$4\MinVecLen_Q(\xi)=\max_{\{v\in V^*\big|[v]=\xi\}} 
Q^*\Big((q(v_0)+2v)\otimes (q(v_0)+2v)\Big)+\rk(V).$$ Since $Q$ is
negative-definite, the maximum exists; and since the number of elements
in the cokernel of $q$ is odd, the
function is independent of the choice of $v_0$.

Let $K$ be an alternating knot.  The projection can be used to give a
planar graph $G$, whose vertices correspond to the white regions in
the checkerboard coloring of the knot projection. In fact, by choosing
a regular alternating projection, we can assume that in the white
graph, there are no edges connecting a vertex to itself. Let $V$
denote the integral lattice formally generated by the vertices, modulo
the relation that $\sum_{v\in V} v=0$.  Recall
from~\cite{BurdeZieschang} that the Goeritz form corresponding to this
projection of $K$ is the quadratic form $$Q\colon V\otimes V
\longrightarrow \Z$$ defined by the rule that if $v$ and $w$ are
distinct vertices of $G$, then $Q(v\otimes w)$ is the number of edges
connecting $v$ to $w$, and also $Q(v\otimes v)=-\deg(v)$, where here
$\deg(v)$ denotes the number of edges containing $v$.  The number of
elements in the cokernel of the associated linear map $q$ is the
determinant $D$ of the knot $K$.

Consider the quadratic form $R_{2n-1}$ on $\Z\oplus \Z$ represented by the matrix
$$\left(\begin{array}{rr}
-n & 1 \\ 1 & -2
\end{array}\right).$$
We define 
$$\phi_0\colon \Zmod{(2n-1)}\stackrel{\cong}\longrightarrow\Coker R_{2n-1}$$
by the property that $\phi_0(i)\in V^*$ is the homomorphism
whose value on $(1,0)$ is $i$ and on $(0,1)$ is $0$.
Correspondingly we define a map
$\gamma_{2n-1}\colon \Zmod{(2n-1)}\longrightarrow \Q$
by
\[
\gamma_{2n-1}(i)=\MinVecLen_{R_{2n-1}}(\phi_0(i)).
\]

Given an alternating projection of a knot $K$ with determinant $D$,
Goeritz form $Q$, an isomorphism $\phi \colon \Zmod{D}\longrightarrow
\Coker(q)$, and a sign $\epsilon=\pm 1$, we define a function
$T_{\phi,\epsilon}\colon \Zmod{D} \longrightarrow \Q$ by $$
T_{\phi,\epsilon}(i)=-\epsilon \cm \MinVecLen_Q(\phi(i))-\gamma_{D}(i).  $$ Note
that such a map satisfies the symmetry
$T_{\phi,\epsilon}(i)=T_{\phi,\epsilon}(-i)$.

\begin{theorem}
\label{thm:AlternatingCriterion}
Let $K$ be an alternating knot with determinant $D$, and let $Q$ be
the negative-definite Goeritz form corresponding to a regular,
alternating projection of $K$.  If $u(K)=1$, then there must be an
isomorphism $\phi\colon
\Zmod{D}\longrightarrow \Coker(q)$ and a sign $\epsilon\in\{\pm 1\}$
with the properties that for all $i\in\Zmod{D}$:
\begin{eqnarray}
T_{\phi,\epsilon}(i)&\equiv& 0 \pmod{2} \label{eq:Integrality}
\\
T_{\phi,\epsilon}(i)&\geq& 0.   \label{eq:Positivity}
\end{eqnarray}
Moreover, if $K$ is an alternating knot with $u(K)=1$,
and furthermore $|\MinVecLen_{Q_K}(0)|\leq \OneHalf$, then there is a choice
of $\epsilon$ and $\phi$ satisfying Equations~\eqref{eq:Integrality}
and~\eqref{eq:Positivity}, and the following additional symmetry:
\begin{equation}
\label{eq:Symmetry}
T_{\phi,\epsilon}(i)=T_{\phi,\epsilon}(2k-i)
\end{equation}
for $1\leq i < k$ when $D=4k-1$ and 
for $0\leq i < k$ when $D=4k+1$.
\end{theorem}

The obstruction given in Theorem~\ref{thm:AlternatingCriterion} does
not depend on the choice of alternating projection of $K$
(c.f. Section~\ref{sec:Proof}).  Note also that there are stronger
versions of the above result,
c.f. Theorem~\ref{thm:SignRefinedVersion} and
\ref{thm:StrongerForm}.

As an application of the above obstruction, we study knots with small
($\leq 10$)
crossing number. Note that if a non-prime knot has crossing number
$\leq 10$, then $H_1(\Sigma(K);\Z)$ is non-cyclic, and hence it does
not have $u(K)=1$.  Alternatively, one can
appeal to a general result of Scharlemann~\cite{Scharlemann},
according to which knots with $u(K)=1$ are
prime.  So, we will always restrict our attention to the prime case here.

We begin with the case of knots with nine or fewer crossings. Among
these, three were listed in Kawauchi's table as having unknown $u=1$
or $2$. For these, we get the following result:

\begin{cor}
\label{cor:VerySmallKnots}
The knots $8_{10}$, $9_{29}$, and $9_{32}$ all have unknotting number
equal to two.
\end{cor}

In recent work~\cite{GordonLueckeUK1}, 
Gordon and Luecke have also shown that $9_{29}$ and $9_{32}$ have $u=2$.

Among knots which admit a diagram with nine or fewer crossings, the
signature, cyclicity of $H_1(\Sigma(K))$, and
Theorem~\ref{thm:AlternatingCriterion} are sufficient to determine all
of those with $u=1$. We review here some earlier results in this
direction. In this family there are 58 knots with $|\sigma|\leq 2$,
c.f.~\cite{Kawauchi} see also~\cite{Kirby}. Of these knots, $35$ can
be directly seen to have $u=1$. Of the remaining 23, eight ($8_{18}$,
$9_{35}$, $9_{37}$, $9_{40}$, $9_{41}$, $9_{46}$, $9_{47}$, and
$9_{48}$) are ruled out by the fact that $H_1(\Sigma(K);\Z)$ is
non-cyclic, see also~\cite{Nakanishi}.  Of the remaining $15$, one
($7_4$) is ruled out by Lickorish~\cite{Lickorish}, nine ($8_3$,
$8_6$, $8_8$, $8_{12}$, $9_5$, $9_8$, $9_{15}$, $9_{17}$, and
$9_{31}$) are ruled out Kanenobu and Murakami~\cite{KanenobuMurakami},
one ($9_{25}$) is ruled out by Kobayashi~\cite{Kobayashi}, and one
more ($8_{16}$) has been ruled out by Rickard, see
also~\cite{MurakamiYasuhara}, \cite{Stoimenow}.  This left open the
unknotting status of $8_{10}$, $9_{32}$, and $9_{29}$ (note that
$9_{29}$ was mistakenly listed in some tables as having $u=1$). While
Corollary~\ref{cor:VerySmallKnots} concerns these last three knots,
same method (Theorem~\ref{thm:AlternatingCriterion}) also proves that
the earlier 12 knots do not have $u=1$.

We also consider knots with ten crossings\footnote{In this
paper, we
use the numbering scheme on knots from Rolfsen's table~\cite{Rolfsen},
modified so that the Perko pair is removed. 
In particular, $10_{161}\neq 10_{162}$, and
the pair  $10_{83}$ and $10_{86}$ have been switched by comparison with Kawauchi's table.}.  There are many knots in
Kawauchi's table~\cite{Kawauchi} with unknown unknotting status. We
considered the knots in the table which were listed as having
unknotting number undetermined but possibly one. There are 40 such
knots, of which 28 are alternating and 12 are
not. Theorem~\ref{thm:AlternatingCriterion} shows that none of these
$28$ alternating knots has $u(K)=1$:

\begin{cor}
\label{cor:Alt10}
The following twenty-four
ten-crossing alternating knots have unknotting number equal to
two:
$$\begin{array}{lllllllllllll}
10_{48}  & 10_{52} & 10_{57} & 10_{58} & 10_{64} &
10_{67} & 10_{68} & 10_{70} & 10_{81} & 10_{83} &
10_{86} & 10_{87} \\
10_{90} &
10_{93} & 10_{94} &
10_{96} &  10_{105} & 10_{106} & 10_{109} &
10_{110} &  10_{112} & 10_{116} & 10_{117} &
10_{121};
\end{array}
$$
moreover, the knots
$$
\begin{array}{llll}
10_{51} & 10_{54} & 10_{77} & 10_{79}
\end{array}
$$
have unknotting number equal to two or three.
\end{cor}

Note that in an earlier  preprint using the linking form, 
Stoimenow~\cite{Stoimenow} has 
shown that the knots 
$$ 10_{86}, 10_{105}, 10_{106}, 10_{109}, 10_{116}, 10_{121}. $$
have $u=2$. In recent work, Gordon and Luecke~\cite{GordonLueckeUK1}
have shown the same results for
\[10_{79}, 10_{81}, 10_{83}, 10_{86}, 10_{87}, 10_{90}, 10_{93}, 10_{94}, 10_{96}.\]

Theorem~\ref{thm:AlternatingCriterion}, and the classical
invariants (signature, $H_1(\Sigma(K);\Z)$) suffice to classify all
$10$-crossing alternating knots with $u(K)=1$,
see also Section~\ref{sec:Calcs}.

Of the remaining twelve non-alternating ten-crossing knots in
Kawauchi's table with unknown unknotting number equal possibly to one,
one ($10_{145}$) has been shown to have $u=2$ by Tanaka~\cite{Tanaka}.
Note that $u(10_{131})=1$ (c.f.~\cite{Stoimenow}, see also
Figure~\ref{fig:10s131} below).  We resolve here the status of an
additional 9 of these, using the methods for the proof of
Theorem~\ref{thm:AlternatingCriterion}, which apply for certain
non-alternating cases, c.f. Section~\ref{sec:NonAlternating} below.

\begin{cor}
\label{cor:Rest10}
The following nine non-alternating, ten-crossing knots 
$$
\begin{array}{cccccccccc}
10_{125} & 10_{126} & 10_{130} & 10_{135} & 10_{138} &
10_{148} & 10_{151} & 10_{158} & 10_{162}
\end{array}
$$
have unknotting number equal to two.
\end{cor}

Recent work of Gordon and Lucke~\cite{GordonLueckeUK1} shows that
\[ 10_{148},\hskip1cm 10_{151},\hskip1cm 10_{153} \]
have $u=2$. Thus, the results of this paper together with~\cite{GordonLueckeUK1} completes the classification of all $10$-crossing
knots with $u=1$.

\subsection{The basic idea}
We discuss some of the ingredients which go into the proof of
Theorem~\ref{thm:AlternatingCriterion}.
We begin with the following observation of Montesinos
c.f.~\cite{Montesinos}, \cite{Lickorish}, \cite{Bleiler}:

\begin{lemma} (Montesinos)
\label{lemma:Montesinos}
If $K$ has $u(K)=1$, then
$\Sigma(K)\cong S^3_{\pm D/2}(C)$ for some other knot $C\subset S^3$,
where here $D$ is the determinant of $K$.
\end{lemma}

In using the above lemma, it is helpful to have an invariant for
three-manifolds which one can calculate for a three-manifold given as
a branched double-cover of $S^3$ branched along a specific knot, and
which can also detect obstructions to realizing a given three-manifold
$Y$ as $n/2$ surgery on a knot in $S^3$ (for integral $n$). Such an
obstruction is furnished by Heegaard Floer homology.

The algebraic structure of Heegaard Floer homology, together with
an induced grading (which takes values in the rational numbers $\Q$)
gives rise to a function
$$d\colon H^2(\Sigma(K);\Z)\longrightarrow \Q,$$
the ``correction terms''
for $\Sigma(K)$ (c.f.~\cite{AbsGraded}, see also the discussion in Section~\ref{sec:Review}
below). (Indeed, this map can be
given for an arbitrary oriented rational homology-three-sphere $Y$
instead of $\Sigma(K)$, except that in the general case, the
correction terms should be interpreted as a rational-valued function
on the $\SpinC$ structures over $Y$.)
The correction terms constrain the intersection form of any smooth
four-manifold which bounds $\Sigma(K)$, according
to Theorem~\ref{AbsGraded:thm:IntFormQSphere} of~\cite{AbsGraded}
(restated in Theorem~\ref{thm:IntForms} below), which is analogous to 
a gauge-theoretic result of Fr{\o}yshov~\cite{Froyshov}.

This fact, together with Lemma~\ref{lemma:Montesinos}, leads at once
to an obstruction to a knot having unknotting number one, stated in
terms of the correction terms of $\Sigma(K)$, c.f.
Theorem~\ref{thm:IntFormBounds} below. Note that that result does not
use the hypothesis that $K$ is alternating.

In general, calculating Heegaard Floer homology, and even the
correction terms $d(Y,\spinc)$ for an arbitrary three-manifold can be
quite challenging. However, this problem is easily solved for the
branched double-covers of alternating knots, and in particular the
correction terms correspond to the quantities $\MinVecLen_Q$ for the
Goeritz form (see Proposition~\ref{BrDCov:prop:AltLink} and
Theorem~\ref{BrDCov:thm:AbsGradeAlt} both in ~\cite{BrDCov}, restated
as Proposition~\ref{prop:AltKnots} below).

In fact, the calculation of the correction terms for alternating
knots, together with the intersection form bounds are not
sufficient for establishing the full statement of
Theorem~\ref{thm:AlternatingCriterion}: they establish
only Conditions~\eqref{eq:Integrality} and \eqref{eq:Positivity}. For 
the additional symmetry from Condition~\eqref{eq:Symmetry},
we need to go  further into the structure of the Floer
homology of branched double-covers of alternating knots.
Specifically, these are rational homology three-spheres whose Heegaard
Floer homology is as simple as possible: they are $L$-spaces in the
sense of~\cite{NoteLens} and also Definition~\ref{def:LSpace} below.
To obtain the full statement of
Theorem~\ref{thm:AlternatingCriterion}, we establish constraints on
the correction terms of $L$-spaces which can be obtained as $n/2$
surgery on a knot in $S^3$, c.f. Theorem~\ref{thm:LSpaceSymmetry}
below.

\subsection{Comparison with other techniques}

There are other applications of gauge theory and Floer theory to
studying the unknotting number of knots, c.f.~\cite{KMMilnor},
\cite{Rudolph}, \cite{4BallGenus}, \cite{RasmussenThesis}, \cite{OwensStrle}.
These techniques all give various lower bounds on the four-ball genus,
and hence the unknotting number. By contrast, 
the obstructions in this paper give information which is independent
of the four-ball genus, and in particular they also give non-trivial
bounds in cases where the four-ball genus is known to be zero or one.

\subsection{Logical dependence}
The proof of Theorem~\ref{thm:AlternatingCriterion} uses Heegaard
Floer homology, as introduced in~\cite{HolDisk}.  The results here
make heavy use of the surgery long exact sequence for Heegaard Floer
homology, c.f.~\cite{HolDiskTwo}, and also properties of the rational
grading on Floer homology, c.f.~\cite{AbsGraded} (and its interaction
with the long exact sequence). In addition, we do refer some
to~\cite{BrDCov}, a paper concerned with the Heegaard Floer homology
of branched double-covers of knots, and specifically to the results
when $K$ is alternating. But these results are confined 
Section~\ref{BrDCov:sec:AltLink} of~\cite{BrDCov}, and
are all fairly straightforward consequences of the surgery long exact
sequence and the rational gradings. Utilizing in addition some results
from~\cite{NoteLens}, we obtain stronger constraints on alternating
knots with unknotting number one, as explained in 
Theorem~\ref{thm:StrongerForm} below. These stronger
results are not needed, though, for the classification of alternating
knots with unknotting number one with ten or fewer crossings.

\subsection{Organization}
The rest of this paper is organized as follows. In
Section~\ref{sec:Review} we recall some of the essentials of the
Heegaard Floer homology package which we use here, and specifically
the constraints on intersection forms coming from the correction
terms. In Section~\ref{sec:FirstApp}, we show how this leads quickly
to an obstruction for an arbitrary knot having $u(K)=1$, stated in
terms of the correction terms of its branched double-cover (c.f.
Theorem~\ref{thm:IntFormBounds} below). In that section, we also
recall the Heegaard Floer homology of branched double-covers of
alternating knots, including the interpretation of their correction
terms in terms of the Goeritz matrix.  In
Section~\ref{sec:LSpaceSurgeries}, we prove a result about the
correction terms of an $L$-space which can be obtained as $D/2$
surgery on a knot in $S^3$ (for some integer $D$). This provides a key
ingredient in establishing Condition~\eqref{eq:Symmetry} from
Theorem~\ref{thm:AlternatingCriterion}. Theorem~\ref{thm:AlternatingCriterion}
is then proved in Section~\ref{sec:Proof}.  In Section~\ref{sec:Calcs}
we discuss the calculations which lead to the proofs of
Corollaries~\ref{cor:VerySmallKnots} and \ref{cor:Alt10}.  In
Section~\ref{sec:NonAlternating}, we turn our attention to some
non-alternating knots, proving Corollary~\ref{cor:Rest10} above. In
Section~\ref{sec:Refine} we give refinements of
Theorem~\ref{thm:AlternatingCriterion}. 
In particular, we show how the methods shed light on the problem of signed
crossing (which we illustrate with the knot $9_{33}$).
Moreover, in that section we
give an interpretation of the non-negative integers
$T_{\phi,\epsilon}$.

\subsection{Acknowledgements}

The authors wish to thank Cameron Gordon, Akio Kawauchi, 
Charles Livingston, Jacob
Rasmussen, and Adam Sikora for helpful conversations and correspondence.

\section{Heegaard Floer homology}
\label{sec:Review}

We give here a rapid outline of the Heegaard Floer homology needed in
the present article.  We consider oriented three-manifolds $Y$ which are
rational homology three-spheres (i.e. closed three-manifolds
with $H_1(Y;\Q)=0$), and for simplicity, we use here Heegaard Floer
homology with coefficients in a field $\Field$, which we take to be
$\Zmod{2}$ for definiteness.

\subsection{Heegaard Floer homology for rational homology three-spheres
and its $\Q$-grading}

Recall that if $X$ is an oriented three- or four-manifold, the space
of $\SpinC$ structures $\SpinC(X)$ is an affine space for $H^2(X;\Z)$.
Moreover, each $\SpinC$ structure has a first Chern class in
$H^2(X;\Z)$, which is related to the action by the formula
$c_1(\spinc+h)=c_1(\spinc)+2h$ for any $h\in H^2(X;\Z)$.

When $Y$ is an oriented rational homology three-sphere and $\spinc$ is
a $\SpinC$ structure over $Y$, its Heegaard Floer homology
$\HFp(Y,\spinc)$ is a $\Q$-graded module over the polynomial algebra
$\Field[U]$, 
$$\HFp(Y,\spinc)=\bigoplus_{d\in\Q} \HFp_d(Y,\spinc),$$
where multiplication by $U$ lowers degree by two. 
In each grading, $i\in\Q$, $\HFp_i(Y,\spinc)$ is a 
finite-dimensional $\Field$-vector space.

Indeed, there is another simpler variant of Heegaard Floer homology,
$\HFinf(Y)$, for which
\begin{equation}
\label{eq:HFinf}
\HFinf(Y,\spinc)\cong\Field[U,U^{-1}]
\end{equation}
for each $\spinc\in\SpinC(Y)$
(c.f. Theorem~\ref{HolDiskTwo:thm:HFinfGen} of~\cite{HolDiskTwo}), and
which admits a natural $\Field[U]$-equivariant map $$\pi \colon
\HFinf(Y,\spinc)\longrightarrow \HFp(Y,\spinc)$$ which is zero in all
sufficiently negative degrees and an isomorphism in all sufficiently
positive degrees. Note that the quotient
$$\HFpRed(Y)=\HFp(Y)/\pi(\HFinf(Y))$$ is a finite-dimensional
$\Field$-vector space.

The image of $\pi$ determines a function $$d\colon \SpinC(Y)
\longrightarrow \Q$$ (the ``correction terms'' of~\cite{AbsGraded}) which associates to each $\SpinC$ structure the
minimal $\Q$-grading of any (non-zero) homogeneous element in
$\HFp(Y,\spinc)$ in the image of $\pi$.  Note that orientations play a
vital role in Heegaard Floer homology. For example, the correction
terms for $Y$ and its opposite $-Y$ are related by the formula
\begin{equation}
\label{eq:ReverseOrientation}
d(-Y,\spinc)=-d(Y,\spinc),
\end{equation} under a natural identification
$\SpinC(Y)\cong \SpinC(-Y)$.

There is a conjugation symmetry on the space of $\SpinC$ structures
$\spinc\mapsto {\overline\spinc}$. Heegaard Floer homology is invariant
under this symmetry, in the sense that we have a commutative square
$$\begin{CD}
\HFinf(Y,\spinc)@>{\cong}>>\HFinf(Y,{\overline\spinc}) \\
@V{\pi}VV @V{\pi}VV \\
\HFp(Y,\spinc)@>{\cong}>> \HFp(Y,{\overline\spinc}) \\
\end{CD}
$$
In particular, we have that
\begin{equation}
\label{eq:ConjugationSymmetry}
d(Y,\spinc)=d(Y,{\overline \spinc}).
\end{equation}

\subsection{$\Zmod{2}$ gradings}
It is sometimes convenient to consider $\HFp(Y,\spinct)$ as a
$\Zmod{2}$-graded (rather than $\Q$-graded) theory. 
Let $X$ be a smooth, oriented four-manifold with boundary
diffeomorphic to $Y$, equipped with a $\SpinC$ structure $\spinc$ with
$\spinc|_{Y}=\spinct$.  The quantity
$$\frac{c_1(\spinc)^2-\sigma(X)}{4},$$
thought of as an element in
$\Q/2\Z$ is easily seen to be independent of the extending
four-manifold $X$. If $\xi\in \HFp_i(Y)$ is a non-zero element
for some $i\in\Q$, 
then the difference
$$i-\frac{c_1(\spinc)^2-\sigma(X)}{4}$$ is an integer 
(for any
choice of $X$), and indeed its parity 
(which is independent of $X$) is the $\Zmod{2}$ grading of
$\xi$.

\subsection{Naturality under cobordisms}

Let $X$ be a smooth, connected,
oriented four-manifold with boundary given by $\partial X=-Y_0\cup
Y_1$ where $Y_0$ and $Y_1$ are connected, oriented three-manifolds. We
call such a four-manifold a cobordism from $Y_0$ to $Y_1$.  If $X$ is
a cobordism from $Y_0$ to $Y_1$, and $\spinc\in\SpinC(X)$ is a $\SpinC$
structure, then there are naturally induced maps on Heegaard Floer homology
which fit into the following diagram:
$$\begin{CD}
\HFinf(Y_0,\spinc_0)@>{F^{\infty}_{X,\spinc}}>> \HFinf(Y_1,\spinc_1) \\
@V{\pi_0}VV @VV{\pi_1}V \\
\HFp(Y_0,\spinc_0)@>{F^{+}_{X,\spinc}}>> \HFp(Y_1,\spinc_1), \\
\end{CD}$$
where here $\spinc_i$ denotes the restriction of $\spinc$ to $Y_i$.
For fixed $X$ and $\xi\in\HFp(Y_0)$, we have that
$\HFp_{X,\spinc}(\xi)=0$ for all but finitely many
$\spinc\in\SpinC(X)$, c.f. Theorem~\ref{HolDiskFour:thm:Finiteness}
of~\cite{HolDiskFour}, and hence there is a well-defined map
$$F^+_{X}\colon \HFp(Y_0)
\longrightarrow \HFp(Y_1)$$ defined by
$$F^+_{X}=\sum_{\spinc\in\SpinC(X)}F^+_{X,\spinc}$$ (note that the
same construction does not work for $\HFinf$: for a given
$\xi\in\HFinf(Y_0)$, there might be infinitely many different
$\spinc\in\SpinC(X)$ for which $F^\infty_{X,\spinc}(\xi)$ is
non-zero).

The map $\pi\colon \HFinf(Y,\spinc) \longrightarrow \HFp(Y,\spinc)$
preserves the $\Q$-grading, and moreover, maps induced
by cobordisms $F^\circ_{X,\spinc}=F^{\infty}_{X,\spinc}$ or
$F^+_{X,\spinc}$ respect the $\Q$-grading in the following sense. If
$Y_0$ and $Y_1$ are rational homology three-spheres, and 
$X$ is a cobordism from $Y_0$ to $Y_1$, with $\SpinC$ structure $\spinc$,
the map induced by the cobordism maps
$$F^{\circ}_{X,\spinc}\colon \HFc_d(Y_0,\spinc_0) \longrightarrow
\HFc_{d+\Delta}(Y_1,\spinc_1),$$
for
\begin{equation}
\label{eq:GradingShift}
\Delta = \frac{c_1(\spinc)^2-2\chi(X)-3\sigma(X)}{4},
\end{equation}
where here $\HFc=\HFinf$ or $\HFp$,
$\chi(X)$ denotes the Euler characteristic of $X$, and
$\sigma(X)$ denotes its signature.  In fact
(c.f. Theorem~\ref{HolDiskFour:thm:AbsGrade} of~\cite{HolDiskFour})
the $\Q$ grading is uniquely characterized by the above property,
together with the fact that $d(S^3)=0$.

Naturality of the maps induced by cobordisms can be phrased as follows.
Suppose that $W_0$ is a smooth cobordism from $Y_0$ to $Y_1$ and $W_1$
is a cobordism from $Y_1$ to $Y_2$, and suppose moreover that $Y_i$ are
rational homology three-spheres, then 
$$\Fp{W_1}\circ \Fp{W_0} = \Fp{W_0\cup_{Y_1} W_1}$$
(c.f. Theorem~\ref{HolDiskFour:thm:Composition} of~\cite{HolDiskFour}).

\subsection{Intersection form bounds}

The correction terms of a rational homology three-sphere $Y$ constrain
the intersection forms of smooth four-manifolds which bound $Y$,
according to the following result, which is analogous to a gauge-theoretic
result of Fr{\o}yshov~\cite{Froyshov}:

\begin{theorem}
  \label{thm:IntForms} Let $Y$ be a rational homology and $W$ be a
  smooth four-manifold which bounds $Y$ with negative-definite
  intersection form. Then, for each $\SpinC$ structure $\spinc$
  over $W$, we have that 
\begin{equation} 
        \label{eq:CongIF}
        c_1(\spinc)^2+b_2(W)\equiv
  4d(Y,\spinc)\pmod{2}, 
\end{equation} 
and we have the inequality
\begin{equation}
        \label{eq:Pos}
        c_1(\spinc)^2+b_2(W)\leq
        4 d(Y,\spinc|_{Y}).
\end{equation}
\end{theorem}

The proof of the above theorem can be found in
Theorem~\ref{AbsGraded:thm:IntFormQSphere} of~\cite{AbsGraded}. In the
case where $b_1(W)=0$, Equation~\eqref{eq:CongIF} follows easily from
the characterization of the $\Q$ grading
(c.f. Equation~\eqref{eq:GradingShift} above), while
Inequality~\eqref{eq:Pos} follows quickly from the fact (c.f. the
proof of Theorem~\ref{AbsGraded:thm:Donaldson} in~\cite{AbsGraded})
that if $X$ is a cobordism with $b_2^+(X)=b_1(X)=0$ between rational
homology three-spheres $Y_0$ to $Y_1$ and $\spinc\in\SpinC(X)$ then
the induced map $F^{\infty}_{X,\spinc}$ is an isomorphism.  In the
case where $b_1(W)>0$, one can perform surgery to reduce to the
previous case.

By contrast, if $X$ is a cobordism from $Y_0$ to $Y_1$ with $b_2^+\neq
0$, then the induced map on $\HFinf$ is trivial,
c.f. Lemma~\ref{HolDiskFour:lemma:BTwoPlusLemma}
of~\cite{HolDiskFour}.

\subsection{Long exact sequences}

Heegaard Floer homology satisfies a surgery long exact sequence, which we state presently.
Suppose that $M$ is a three-manifold with torus boundary, and fix
three simple, closed curves $\gamma_0$, $\gamma_1$, and $\gamma_2$ in
$\partial M$ with 
\begin{equation}
\label{eq:TriadRelation}
\#(\gamma_0\cap \gamma_1)=\#(\gamma_1 \cap \gamma_2)
=\#(\gamma_2\cap \gamma_0) = -1
\end{equation} (where here the algebraic intersection
number is calculated in $\partial M$, oriented as the boundary of
$M$), so that $Y_0$ resp. $Y_1$ resp. $Y_2$ are obtained from $M$ by
attaching a solid torus along the boundary with meridian $\gamma_0$
resp. $\gamma_1$ resp. $\gamma_2$.
Note that there are two-handle cobordisms $W_i$ connecting $Y_i$ to $Y_{i+1}$
(where we view $i\in\Zmod{3}$).

\begin{theorem}
\label{thm:ExactSeq}
Let $Y_0$, $Y_1$, and $Y_2$ be related as above. Then, we have a long exact sequence
$$
\begin{CD}
...@>>>\HFp(Y_0) @>>>\HFp(Y_1) @>>>\HFp(Y_2)@>>>...
\end{CD}
$$
where the maps are induced by from the natural two-handle cobordisms.
\end{theorem}

The above theorem is proved in Theorem~\ref{HolDiskTwo:thm:GeneralSurgery} of~\cite{HolDiskTwo}.

\subsection{The case where $b_1(Y)>0$}

Although we have restricted our attention mainly to the case of
rational homology three-spheres, the construction of Heegaard Floer
homology (and in particular Theorem~\ref{thm:ExactSeq}) works for
arbitrary closed, oriented three-manifolds. The rational grading which
we discussed here, however, works only for $\SpinC$ structures
whose first Chern class is torsion. Also, the
structure of $\HFinf(Y,\spinc)$ is slightly more intricate. We state
now the case where $H_1(Y;\Z)\cong\Z$. In this case, the map $\pi$ is
trivial for all $\SpinC$ structures with non-zero first Chern class,
and in the case where $c_1(\spinc)=0$, we have that
\begin{equation}
\label{eq:HFinfBOne}
\HFinf(Y,\spinc)\cong\Field[U,U^{-1}]\oplus \Field[U,U^{-1}];
\end{equation}
the first summand has $\Q$-grading in $\OneHalf + 2\Z$, while the second
has $\Q$-grading in $-\OneHalf+2\Z$.

\subsection{$L$-spaces}

There is a class of three-manifolds for which the Floer homology
is particularly simple.

\begin{defn}
\label{def:LSpace}
A rational homology
three-sphere is called an $L$-space if for each $\spinc\in\SpinC(Y)$,
the map $\pi\colon \HFinf(Y)\longrightarrow \HFp(Y)$ is surjective. 
\end{defn}

Clearly, for an $L$-space, the correction terms determine the Heegaard
Floer homology $\HFp(Y)$. The reader can find another equivalent
definition in~\cite{NoteLens}, which uses a different variant of
Heegaard Floer homology. Finally, it should be pointed out that Floer
homology depends on the choice of coefficient system, so it would be
more precise to call a three-manifold an $L$-space with coefficients
in $\Field$ (our underlying coefficient system), but we do not do this
here.

Note that if $Y$ is an $L$-space, then so is $-Y$.

The following principle gives a plentiful supply of
$L$-spaces. 

\begin{prop}
\label{prop:LSpaces}
Suppose that $Y_0$, $Y_1$, and $Y_2$ are three rational
homology three-spheres related by an exact sequence as in
Theorem~\ref{thm:ExactSeq}, and suppose that the number of elements in
$H_1(Y_1;\Z)$ is greater than the number of elements in $H_1(Y_0;\Z)$
and $H_1(Y_2;\Z)$. Suppose moreover that $Y_0$ and $Y_2$ are
$L$-spaces, then so is $Y_1$.
\end{prop}

The above proposition is straightforward application of
Theorem~\ref{thm:ExactSeq}, together with the structure of $\HFinf$
for a rational homology three-sphere (Equation~\eqref{eq:HFinf}).
Details can be found in Proposition~\ref{NoteLens:prop:LSpaces} of~\cite{NoteLens}.

\subsection{Sharp four-manifolds}
\label{subsec:Sharp}

We will describe here a mechanism which is sometimes useful in the
calculation of correction terms for $L$-spaces. The methods employed
here are closely related to discussions in~\cite{SomePlumbs}
and also Section~\ref{BrDCov:sec:AltLink} of~\cite{BrDCov}.

\begin{defn}
Let $Y$ be an $L$-space. A smooth four-manifold
$X$ which bounds $-Y$ is called {\em sharp} if $b_1(X)=0$,
the intersection form on $H^2(X,Y;\Z)$ is
negative-definite, and if for each $\spinct\in\SpinC(Y)$, there
is some extension $\spinc$ over $X$ with the property that
\begin{equation} 
        \label{eq:Cong}
        c_1(\spinc)^2+b_2(X)=-4d(Y,\spinct).
\end{equation}
\end{defn}

Clearly, if $X$ is sharp, then so is its blowup $X\#\mCP$.

\begin{prop}
\label{prop:Sharp}
Suppose that $Y_0$, $Y_1$, and $Y_2$ are three rational homology
three-spheres related by an exact sequence as in
Theorem~\ref{thm:ExactSeq}, and that $Y_0$ and $Y_2$ are $L$-spaces.
Suppose that there is a four-manifold $X_2$ which is sharp for $Y_2$;
and also that $X_0=W_0\cup W_1\cup X_2$ is sharp for $Y_0$. Then $X_1=W_1\cup
X_2$ is sharp for $Y_1$. (Here, $W_i$ are the two-handle cobordisms
connecting $Y_i$ to $Y_{i+1}$.)
\end{prop}

Suppose that $X$ is negative-definite four-manifold with $b_1(X)=0$
which bounds an $L$-space $-Y$.  Then $X-B^4$ gives a cobordism from
$Y$ to $S^3$. Define $\FormIFp(Y;X)$ to be the set of maps
$$\phi\colon \SpinC(X) \longrightarrow \HFp(S^3)$$ satisfying the
relation that
\begin{equation}
\label{eq:DefiningRelation}
 U^{\frac{\langle c_1(\spinc),v\rangle + v\cm v}{2}} \cm \phi(\spinc+\PD(v))=\phi(\spinc),
\end{equation}
for any $\spinc\in\SpinC(X)$ and $v\in H_2(X;\Z)$,
for which $\langle c_1(\spinc),v\rangle + v\cm v \geq 0$.
The $\Field[U]$-module structure on $\HFp(S^3)$
gives this set  the structure of a $\Field[U]$ module.

There is a map $$T^+_X\colon \HFp(Y) \longrightarrow
\FormIFp(Y;X)$$
given by $$\langle T^+_X(\xi),[\spinc]\rangle
= \Fp{X-B^4,\spinc}(\xi),$$ thinking of $X-B^4$ as a cobordism from
$Y$ to $S^3$. 
We claim that the image of $T^+_X$ is in 
$\FormIFp(Y;X)\subset \Hom(\SpinC(X),\HFp(S^3))$. This follows at once from the fact that a
negative-definite cobordism with $b_1=0$
between a rational homology spheres
induces an isomorphism on $\HFinf$ for each $\SpinC$ structures, whose
degree is given by Equation~\eqref{eq:GradingShift}.

Given $d\in\Q$, let $\FormIFp_d(Y;X)\subset \FormIFp(Y;X)$
denote the set of maps $\phi$ the property that for all
$\spinc\in\SpinC(X)$ for which $\phi(\spinc)\neq 0$, we have that
$\phi(\spinc)\in\HFp_{i}(S^3)\subset \HFp(S^3)$ where
$$i-\left(\frac{c_1(\spinc)^2+2\cm\chi(X)+3\cm \sigma(X)}{4}\right)=d.$$
Clearly, the restriction of $T^+_{X}$ to $\HFp_d(Y)$ is contained in
$\FormIFp_d(Y;X)$.

Clearly, 
$X$ is a sharp four-manifold if and only if the map $T^+_X$ is an
isomorphism.

\begin{lemma}
\label{lemma:AlgStructure}
Let $Y$ be a rational homology three-sphere, and 
suppose that $X$ is a smooth four-manifold with
negative-definite intersection form.
Then,
there is an isomorphism
$$\FormIFp(Y;X)\cong \left(\Field[U,U^{-1}]/\Field[U]\right)^{n}, $$
where here $n$ is the number of elements in the image
of $H^2(X;\Z)$ inside $H^2(Y;\Z)$, and the isomorphism is to be viewed
as an isomorphism between ungraded $\Field[U]$-modules.
\end{lemma}

\begin{proof}
  The isomorphism is induced as follows. Choose for each
  $\spinc\in\SpinC(Y)$ which extends over $X$ an extension
  ${\widetilde \spinc}\in\SpinC(X)$ for which $c_1({\widetilde
    \spinc})^2$ is maximal. This can be done since $X$ has
  negative-definite intersection form. The map from $\FormIFp(Y;X)$
  maps $\phi$ to $(\phi({\widetilde \spinc}_1),...,\phi({\widetilde
    \spinc}_n))$, where $\{\spinc_i\}_{i=1}^n$ is the set of
  $\SpinC$ structures over $Y$ which extend over $X$. It is straightforward
  to write down an inverse for this map.
\end{proof}

Let $Y_0$ and $Y_1$ be $L$-spaces. 
Let $W_0$ is a negative-definite cobordism with $b_1(W_0)=0$ 
from $Y_0$ to $Y_1$,
and let $X_1$ be a negative-definite four-manifold with
$b_1(X_1)=0$ which
bounds $-Y_1$.
Then, there is  an induced map:
$$
\FormF{W_0}\colon \FormIFp(Y_0;W_0\cup X_1) \longrightarrow \FormIFp(Y_1;X_1).
$$
defined by
$$\langle \FormF{W_0}(\phi),\spinc_1\rangle = 
\sum_{\{\spinc\in\SpinC(W_0\cup_{Y_1} X_1)\big| \spinc|_{X_1}=\spinc_1\}}
\phi(\spinc).$$
This is easily seen to be a
well-defined since $b_2^+(W_0\cup X_1)=0$ and $\phi$ satisfies 
Equation~\eqref{eq:DefiningRelation}.
The map is natural under cobordisms, in the sense that the following
square commutes
$$\begin{CD}
\HFp(Y_0) @>{F^+_{W_0}}>>\HFp(Y_1) \\
@V{T^+_{W_0\cup_{Y_1} X_1}}VV @VV{T^+_{X_1}}V \\
\FormIFp(Y_0;W_0\cup_{Y_1} X_1)@>{\FormF{W_0}}>> \FormIFp(Y_1; X_1).
\end{CD}
$$
Commutativity of this square follows at once from the
composition law for the maps induced by cobordisms.

\begin{lemma}
\label{lemma:Inj}
  Suppose that $Y_0$, $Y_1$, and $Y_2$ are three rational homology
  three-spheres related by an exact sequence as in
  Theorem~\ref{thm:ExactSeq}, and suppose that $X_2$ is a
  four-manifold which bounds $-Y_2$ with $b_1(X_2)=0$, 
  so that $W_0\cup W_1\cup X_2$ is a
  negative-definite four-manifold.  Then, the map
  $$\FormF{W_0}\colon \FormIFp(Y_0; W_0\cup W_1\cup
  X_2)\longrightarrow \FormIFp(Y_1;W_1\cup X_2)$$
  is injective.
\end{lemma}

\begin{proof}
  Fix a non-zero $\phi_0\in\FormIFp(Y_0;W_0\cup W_1\cup X_2)$. Its
  non-triviality means that there is an $\spinc\in\SpinC(W_0\cup W_1\cup
  X_2)$ with $\phi_0(\spinc)\neq 0$. Indeed, by multiplying $\phi_0$
  by powers of $U$ if necessary, we obtain a new element
  $\phi\in\FormIFp(Y_0;W_0\cup W_1\cup X_2)$ with the property that
  $U\cm \phi\equiv 0$, but there is some $\spinc\in\SpinC(W_0\cup
  W_1\cup X_2)$ with $\phi(\spinc)\neq 0$. (This follows at once
  from Lemma~\ref{lemma:AlgStructure}.)
  
  The kernel of the map $H^2(W_0\cup W_1\cup X_2;\Z) \longrightarrow
  H^2(W_1\cup X_2;\Z)$ is generated by $\PD(\Sigma_0)$, where
  $\Sigma_0\in H_2(W_0,Y_0;\Z)\cong \Z$ is a generator.  Moreover, the
  kernel of the map $H^2(W_0\cup W_1\cup X_2;\Z)\longrightarrow
  H^2(X_2;\Z)$ is generated by two homology classes $\PD(\Sigma_0)$
  and $\PD(e)$, and $e\in H_2(W_0\cup W_1;\Z)$ with $\Sigma_0\cm
  e=1$, and $e\cm e=-1$. (Indeed, the class $e$ can be
  represented by an embedded two-sphere.)
  
  Fix $\spinc_0$ so that $\phi(\spinc_0)\neq 0$. Indeed,
  by subtracting off $\PD(e)$ if necessary, 
  we can assume without loss of generality that
  $\langle \phi(\spinc_0),e\rangle =-1$.
  Choose a maximal integer $a$ so that
  $\phi(\spinc_0+a\cm \PD(\Sigma_0) + a\cm \PD(e))\neq 0$.
  This exists since $\phi$ has finite support (this in turn follows
  from Equation~\eqref{eq:DefiningRelation}, together with the fact that
  $W_0\cup W_1\cup X_2$ has negative-definite intersection form).
  Note that 
  $\langle c_1(\spinc_0+(a+1)\cdot\PD(\Sigma_0)+a\cdot\PD(e)),e\rangle=+1$, 
  and hence it follows 
  (from Equation~\eqref{eq:DefiningRelation}, together with the choice
  of $a$) that
  \begin{equation}
    \label{eq:ProveMe}
    \phi(\spinc_0+(a+1)\cdot\PD(\Sigma_0)+a\cdot\PD(e))= \phi(\spinc_0+(a+1)\cdot\PD(\Sigma_0)+(a+1)\cdot\PD(e))=0.
  \end{equation}
  Thus, for $\spinc=\spinc_0+a\cm \PD(\Sigma_0)+a\cm \PD(e)$, we have that
  \begin{equation}
    \label{eq:WhatWeWant}
    \langle c_1(\spinc),e\rangle = -1,\hskip1cm
        \phi(\spinc)\neq 0,\hskip1cm \phi(\spinc+\PD(\Sigma_0))=0.
  \end{equation}

  Clearly, in the orbit $\spinc+\Z\cm \PD(\Sigma_0)$, only
  $c_1(\spinc)$ and $c_1(\spinc+\PD(\Sigma_0))$ have evaluation $\pm
  1$ on $e$, and hence only $\spinc$ and $\spinc+\PD(\Sigma_0)$
  have the possibility of having non-trivial
  value under $\phi$. However, Equation~\eqref{eq:WhatWeWant} proves
  that $$\langle \FormF{W_0}(\phi),\spinc|_{W_1\cup X_2}\rangle
  =\sum_{b\in\Z} \phi(\spinc+b\cdot\PD(\Sigma_0))
  =\phi(\spinc)+\phi(\spinc+\PD(\Sigma_0))\neq 0.$$
\end{proof}

\vskip.2cm
\noindent{\bf{Proof of Proposition~\ref{prop:Sharp}.}}
Consider the diagram:
$$\begin{CD}
0@>>>\HFp(Y_0) @>{F_{W_0}}>> \HFp(Y_1) @>{F_{W_1}}>> \HFp(Y_2) @>>> 0 \\
&& @V{\cong}V{T^+_{X_0}}V @VV{T^+_{X_1}}V @V{\cong}V{T^+_{X_2}}V \\
0@>>>\FormIFp(Y_0;X_0) @>{\FormF{W_0}}>> \FormIFp(Y_1;X_1) @>{\FormF{W_1}}>> \FormIFp(Y_2;X_2),
\end{CD}$$
where here the top row is 
exact.
According to Lemma~\ref{lemma:Inj},
$\FormF{W_0}$ is an injection. 
A straightforward diagram-chase now establishes that
$T^+_{X_1}$ is injective.

By Lemma~\ref{lemma:AlgStructure}, $T^+_{X_1}$ must be an isomorphism. This
follows from the following observation:
suppose 
$$f\colon \left(\frac{\Field[U,U^{-1}]}{\Field[U]}\right)^a \longrightarrow \left(\frac{\Field[U,U^{-1}]}{\Field[U]}\right)^b$$
is an injective  map of $\Field[U]$ modules for some $a\geq b$,
then $a=b$ and indeed $f$ is an isomorphism. This can be seen by
restricting to the kernel of $U^n$ (for all $n$),
and appealing to the corresponding fact for finite-dimensional vector spaces.
\qed

\section{First applications of Floer homology}
\label{sec:FirstApp}

In the introduction, there was no reference to $\SpinC$ structures.
To see why these can be eliminated, recall that if
$H^2(Y;\Z)$ has odd order then the map which sends
$\spinc\in\SpinC(Y)$ to half its first Chern class $c_1(\spinc)/2$
induces an isomorphism
$$\SpinC(Y) \longrightarrow H^2(Y;\Z).$$
If $K$ is a knot in $S^3$,
then $H^2(\Sigma(K);\Z)$ has odd order, so we can use the above map to
identify $\SpinC$ structures with integral two-dimensional cohomology
classes. Note that under the above identification, the conjugation symmetry
on $\SpinC(Y)$ is identified with multiplication by $-1$.

With this said, we get the following rather quick consequence of
Lemma~\ref{lemma:Montesinos} and Theorem~\ref{thm:IntForms}, according
to which the correction terms for $\Sigma(K)$ give an
obstruction for $K$ to have unknotting number one.

\begin{theorem}
  \label{thm:IntFormBounds} If $K$ is a knot with unknotting number
  one, then there is some isomorphism $\phi\colon
  \Zmod{D}\stackrel{\cong}{\longrightarrow} H^2(\Sigma(L);\Z)$ with the
  property that at least one of the two conditions holds:
\begin{itemize}
\item   for all $i\in\Zmod{D}$,
\begin{eqnarray*}
\gamma_D(i)
\equiv d(\Sigma(K),\phi(i))\pmod{2} &{\text{and}}&
       \gamma_{D}(i)\leq d(\Sigma(K),\phi(i)),
\end{eqnarray*}
\item   for all $i\in\Zmod{D}$,
\begin{eqnarray*}
\gamma_D(i)\equiv -d(\Sigma(K),\phi(i))\pmod{2} &{\text{and}}&
       \gamma_D(i)\leq -d(\Sigma(K),\phi(i)).
\end{eqnarray*}
\end{itemize}
\end{theorem}

\begin{proof}
If $K$ has unknotting number equal to one, then by Montesinos' trick
(Lemma~\ref{lemma:Montesinos}), $\Sigma(K)\cong \pm S^3_{D/2}(C)$ for
some knot $C\subset S^3$. Thus, $\pm \Sigma(K)$ bounds a four-manifold $W$
with intersection form 
$$\left(\begin{array}{rr}
    -n & 1 \\
    1 & -2 
\end{array}\right),
$$
where here $D=2n-1$.
The theorem now is a direct application of Theorem~\ref{thm:IntForms}
in this context.
\end{proof}

Of course, to apply Theorem~\ref{thm:IntFormBounds} meaningfully, one
must calculate the correction terms for $\Sigma(K)$. When $K$ is an
alternating knot, the correction terms can be calculated in terms of
the Goeritz matrix for $K$. Furthermore, when $K$ is an alternating knot,
$\HFp(\Sigma(K))$ has a particularly simple form.  

\begin{figure}
\mbox{\vbox{\epsfbox{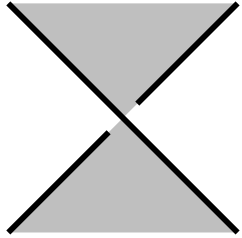}}}
\caption{\label{fig:ColorConventions}
{\bf Coloring conventions for alternating knots.}}
\end{figure}

To this end, we take the Goeritz form, where the knot is colored
according to the coloring conventions as specified in
Figure~\ref{fig:ColorConventions}.

\begin{prop}
\label{prop:AltKnots}
If $K$ is an alternating knot, then $\Sigma(K)$ is an $L$-space;
and indeed, if $Q$ denotes the Goeritz form of $K$ and $q$ the induced map
as described earlier, then
there is an affine identification
$$\phi\colon\Coker (q)\stackrel{\cong}{\longrightarrow} \SpinC(\Sigma(K))
$$
taking zero to the unique spin structure,
with the property that
$$\MinVecLen_Q(\xi)=d(\Sigma(K),\phi(\xi)).$$
\end{prop}

The proof of the above result can be found in~\cite{BrDCov} (combining
Proposition~\ref{BrDCov:prop:AltLink} and
Theorem~\ref{BrDCov:thm:AbsGradeAlt} both in ~\cite{BrDCov}), but
we sketch the argument here for the reader's convenience.

Starting from a connected, alternating projection of a link $L$, one
constructs a four-manifold $X_L$ bounding $\Sigma(L)$ as
follows. For all but one vertex in the black graph, we associate a
one-handle to attach to the four-ball. In terms of Kirby calculus, to
each of these vertices, we associate a dotted unknot. Then, to each
crossing (edge in the black graph), we associate an unknot with
framing $-1$ which links the two dotted unknots which it connects. The
intersection form of $X_L$ is the Goeritz form for $L$, as can be seen
after performing handleslides of the $-1$-framed unknots around the
circuits given by vertices in the white graph. (This description can
be easily seen to agree with the one given in~\cite{BrDCov}, where the
Goeritz form is described with respect to a different basis.)

\begin{figure}
\mbox{\vbox{\epsfbox{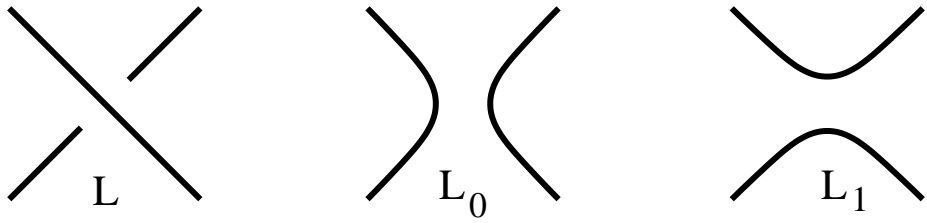}}}
\caption{\label{fig:ResolutionConventions}
{\bf Skein moves.}
Let $L$ be an alternating link, and $L_0$ and $L_1$ be the links obtained
by resolving some fixed crossing of $L$ according to the illustrated
conventions.}
\end{figure}

Fix a crossing for some regular alternating projection of $L$, and let
$L_0$ and $L_1$ be the two links formed by resolving the crossing in
two ways as pictured in Figure~\ref{fig:ResolutionConventions}, and
suppose that both $L_0$ and $L_1$ have connected projection (if such a
point cannot be found, then $L$ is a projection of the unknot).  It is
easy to see that $-\Sigma(L_0)$, $-\Sigma(L_1)$, and $-\Sigma(L)$ are
related as in Theorem~\ref{thm:ExactSeq} (with $Y_0=-\Sigma(L_0)$,
$Y_1=-\Sigma(L)$, and $Y_2=-\Sigma(L_1)$). Now, by induction on the
number of crossings, together with Proposition~\ref{prop:LSpaces},
shows that $\Sigma(L)$ is an $L$-space.  Indeed, observing that
$W_0\cup W_1\cup X_{L_1} = X_{L_0}\# \mCP$, 
the same inductive argument (now using
Proposition~\ref{prop:Sharp}) can be used to show that $X_{L}$ is
sharp for $\Sigma(L)$.

Note that Theorem~\ref{thm:IntFormBounds} and
Proposition~\ref{prop:AltKnots} suffice to establish
Conditions~\eqref{eq:Integrality} and \eqref{eq:Positivity} in
Theorem~\ref{thm:AlternatingCriterion}. However, to establish
Condition~\eqref{eq:Symmetry}, we need an analogous symmetry for the
correction terms of $L$-spaces which are obtained as
$(2n-1)/2$-surgery on a knot in $S^3$.

\section{$L$-space surgeries}
\label{sec:LSpaceSurgeries}

We now turn our attention to the following symmetry result which quickly
provides the missing piece of the proof of 
Theorem~\ref{thm:AlternatingCriterion}.

Note that there is a natural cobordism $W$ from $S^3$ to
$S^3_{-(2n-1)/2}(C)$, whose Kirby calculus picture is given by $C$
with framing $-n$, and a linking unknot with framing $-2$. This
cobordism can naturally be broken into a pair of cobordisms $W_1$ from
$S^3$ to $S^3_{-n}(C)$, followed by a cobordism $W_2$ from
$S^3_{-n}(C)$ to $S^3_{-(2n-1)/2}(C)$. The intersection form for
$W_1\cup_{S^3_{-n}(C)} W_2$ is given by $R_{-(2n-1)/2}$.  Such an
identification is specified by choosing a generator $F$ of
$H_2(W_1;\Z)\cong\Z$, or, equivalently, an orientation on $C$. 
We fix this additional datum for the purposes of the rest of the
present section (though this choice does not affect the final
results, in view of the conjugation symmetry of the correction terms). 

For any knot $C\subset S^3$, we have an isomorphism $$\phi\colon
\Zmod{(2n-1)}\longrightarrow H^2(S^3_{-(2n-1)/2}(C))$$ which takes $i$
to $i\cm \PD[F]|_{S^3_{-(2n-1)/2}(C)}$.  (We suppress the knot $C$
from the notation, but bear in mind that sometimes we use this map
$\phi$ for $C$ and at other times we use $\phi$.)  In the following
statement, we use $c_1(\spinc)/2$ to identify
$\SpinC(S^3_{-(2n-1)/2}(C))\cong H^2(S^3_{-(2n-1)/2}(C);\Z)$.

\begin{theorem}
\label{thm:LSpaceSymmetry}
Let $O$ be the unknot, and 
let $C\subset S^3$ be a knot with the property that for some $n>1$,
$S^3_{-(2n-1)/2}(C)$ is an $L$-space, and suppose moreover that
$$d(S^3_{-(2n-1)/2}(C),0)=d(S^3_{-(2n-1)/2}(O),0).$$
Write $n=2k$ or $2k+1$ for integral $k$.
Then, under the isomorphism $\phi$, we have that
\begin{eqnarray}
\label{eq:dEquality}
\lefteqn{d(S^3_{-(2n-1)/2}(C),\phi(i))-d(S^3_{-(2n-1)/2}(O),\phi(i))} \\
&=&
d(S^3_{-(2n-1)/2}(C),\phi(2k-i))-d(S^3_{-(2n-1)/2}(O),\phi(2k-i))
\nonumber
\end{eqnarray}
for $i=1,...,k$; and when $n=2k+1$, Equation~\eqref{eq:dEquality}
also holds for $i=0$.
\end{theorem}

We give the proof at the end of the present section. 
Indeed, Theorem~\ref{thm:LSpaceSymmetry} follows quickly from another result
(Theorem~\ref{thm:LSpaceFracSurg}) which
identifies the difference in correction terms for the
$-(2n-1)/2$-surgery with a corresponding difference for the
$-n$-surgery. To state this, we introduce some notation.

Let $\{F,S\}$ be the basis for $H_2(W;\Z)$ corresponding to the knot $C$
and the linking unknot (here $S$ is represented by a sphere with square $-2$
and $S\cm F = +1$). 
Let $(a,b)$ be the cohomology class whose evaluation on $F$ and $S$
are $a$ and $b$ respectively.

We now give an ordered list of $2n-1$ characteristic vectors for the
intersection form of $W$, in two cases depending on the parity of $n$.
If $n$ is even, write $n=2k$ and let
\begin{equation}
\label{eq:EvenList}
\kappa_i=
\left\{\begin{array}{ll}
(2i,0) & {\text{for $-k\leq i\leq k$}} \\
(-4k+2i,2) &{\text{for $k<i\leq 2k-1$}} \\
(2i-4k+2,-2) &{\text{for $2k\leq i \leq 3k-2$}} 
\end{array}\right.
\end{equation}
while if $n$ is odd, write $n=2k+1$, and let
\begin{equation}
\label{eq:OddList}
\kappa_i=
\left\{\begin{array}{ll}
(1+2i,-2) & {\text{for $0\leq i \leq k$}} \\
(2i-4k-1,0) & {\text{for $k+1\leq i \leq 3k+1$}} \\
(2i-8k-3,2) &{\text{for $3k+2\leq i \leq 4k$}}
\end{array}\right.
\end{equation}

For $i=0,...,2n-2$, let $w_i\in\SpinC(S^3_{-(2n-1)/2}(C))$ be the
$\SpinC$ structure which can be extended to $\spinc_i\in\SpinC(W)$
with $c_1(\spinc_i)=\kappa_i$.  Let $v_i\in\SpinC(S^3_{-n}(C))$ denote the
restriction of $\spinc_i$ to $S^3_{-n}(C)$.
% Our map from
% $\SpinC(S^3_{-(2n-1)/2}(C))$ to $\SpinC(S^3_{-n}(C))$ sends $w_i$ to
% $v_i$. 

Of course, the $\{v_i\}_{i=0}^{2n-2}$ are not all distinct $\SpinC$
structures on $S^3_{-n}(C)$: the $\SpinC$ structure is determined by
the first coordinate of $\kappa_i\pmod{2n}$.

\begin{theorem}
\label{thm:LSpaceFracSurg}
Assume that $C\subset S^3$ is a knot with the property that for some $n>1$
$S^3_{-(2n-1)/2}(C)$ is an $L$-space. Then so is $S^3_{-n}(C)$.
Moreover, if 
$$d(S^3_{-(2n-1)/2}(C),w_0)=d(S^3_{-(2n-1)/2}(O),w_0),$$
then we have for all $i=0,...,2n-2$,
$$d(S^3_{-n}(C),v_i) - d(S^3_{-n}(O),v_i)
= d(S^3_{-(2n-1)/2}(C),w_i)-d(S^3_{-(2n-1)/2}(O),w_i).$$
\end{theorem}

The proof will be given in the end of the section.

\begin{lemma}
\label{lemma:Positivity}
Let $C$ be any knot in $S^3$.  For $i\in\Zmod{(2n-1)}$,
\begin{equation}
\label{eq:dIneq}
0\leq d(S^3_{-n}(C),v_i)-d(S^3_{-n}(O),v_i)
\leq d(S^3_{-(2n-1)/2}(C),w_i) -d(S^3_{-(2n-1)/2}(O),w_i) .
\end{equation}
\end{lemma}

\begin{proof}
The argument which establishes
Theorem~\ref{thm:IntForms} actually proves the following result:
suppose that $W$ is a cobordism from $Y_0$ to $Y_1$, both of which are
rational homology three-spheres, and suppose that the intersection
form of $W$ is negative-definite, then for each $\spinc\in\SpinC(W)$,
\begin{equation}
\label{eq:MonotoneD}
d(Y_0,\spinc|_{Y_0}) + \left(\frac{c_1^2(\spinc)-2\chi(W)-3\sigma(W)}{4}\right)
\leq d(Y_1,\spinc|_{Y_1}).
\end{equation}
(This can alternatively thought of as a consequence of
Theorem~\ref{thm:IntForms}, together with the additivity of the
correction terms under connected sums, c.f. 
Theorem~\ref{AbsGraded:thm:AdditivityOfD}
of~\cite{AbsGraded}.)

This establishes the existence of constants $k(n,i)$ with the property
that for any knot $C$, $$d(S^3_{-n}(C),w_i) + k(n,i) \leq
d(S^3_{-(2n-1)/2}(C),v_i).$$ Indeed, we can choose the constants so
that equality holds when $C=O$. This follows readily from the fact
that the negative-definite plumbing for $S^3_{-(2n-1)/2}(O)$ is sharp,
which in turn is an easy
consequence of Proposition~\ref{prop:Sharp}, see
also~\cite{SomePlumbs}.  This establishes the second inequality.  The
first is proved along the same lines.
\end{proof}

The following result about $L$-spaces will be useful to us, as well. A
more general result is established in Section~\ref{KMOS:sec:Rational}
of~\cite{KMOS} in the context of the Seiberg-Witten monopole Floer homology
(though the proof adapts readily to the case at hand).

\begin{prop}
\label{prop:LargerLSpaces}
Suppose that $S^3_{-(2n-1)/2}(C)$ is an $L$-space. Then so is $S^3_{-n}(C)$.
\end{prop}

\begin{lemma}
\label{lemma:Surjectivity}
If $C$ is any knot in $S^3$, then the map induced by the cobordism
$$
\begin{CD}
\HFp(S^3_{(2n-1)/2}(C))@>>>\HFp(S^3_{n}(C)) 
\end{CD}
$$ is surjective.
\end{lemma}

\begin{proof}
First, we argue that for all $n\geq 2$, the natural map 
$\HFp(S^3_{n-1}(C)) \longrightarrow \HFp(S^3_{n}(C))$ is injective for all
$n$. This follows at once from the surgery exact sequence
$$\begin{CD}
...@>{D}>>\HFp(S^3_{n-1}(C)) @>{A}>>\HFp(S^3_n(C)) @>>>\HFp(S^3)@>{D}>>...,
\end{CD}$$
together with the observation that $D$ is induced by a cobordism $W$
with $b_2^+(W)=1$, and hence the induced map from $\HFp(S^3)$ is trivial.
Thus, we see that the surgery long exact sequence becomes a short exact
sequence.

Next, consider the long exact sequence
$$
\begin{CD}
@>{B}>>\HFp(S^3_{n-1}(C)) @>>>\HFp(S^3_{(2n-1)/2}(C)) @>>>\HFp(S^3_n(C))@>{B}>>
\end{CD}
$$
Observe that the cobordism $A\circ B$ admits an alternative
factorization as the standard (two-handle) cobordism from $S^3_{n}(C)$
to $S^3_{n}(C)\#(S^2\times S^1)$, followed by another cobordism (in
which the generator of $S^2\times S^1$ becomes null-homologous). From
this, it follows readily that the induced $A\circ B=0$ 
(c.f. Lemma~\ref{KT:lemma:CompositeLemma} of~\cite{calcKT}). Since $A$ is injective, it follows
that $B=0$. Exactness now shows that the stated map is surjective.
\end{proof}

\vskip.2cm
{\noindent{\bf Proof of Proposition~\ref{prop:LargerLSpaces}.}}
Lemma~\ref{lemma:Surjectivity} proves at once that 
if $S^3_{2n-1/2}(C_1)$ is an $L$-space, then so is $S^3_{n}(C_1)$.
Note that $Y$ is an $L$-space iff $-Y$ is, and also if
$C_1$ is the mirror of $C$,  then
$S^3_{r}(C_1)=-S^3_{-r}(C)$ for any rational number $r$.
Thus, the claim follows.
\qed
\vskip.2cm

If $Y$ is a rational homology three-sphere, write
$$d(Y) = \sum_{\spinc\in\SpinC(Y)} d(Y,\spinc).$$

\begin{prop}
\label{prop:OurEuler}
Let $C\subset S^3$ be a knot, and suppose that $S^3_{-(2n-1)/2}$
is an $L$-space. Then, 
$$d(S^3_{-(2n-1)/2}(C))-2 \cdot d(S^3_{-n}(C))
= 
d(S^3_{-(2n-1)/2}(O))-2 \cdot d(S^3_{-n}(O)).$$
\end{prop}

The above follows from a more general result (Proposition~\ref{prop:Euler})
proved in Subsection~\ref{subsec:Euler}.

With the preliminaries in place, we can now turn to the proofs of the
two theorems stated in the beginning of the present section. 

\vskip.2cm
\noindent{\bf{Proof of Theorem~\ref{thm:LSpaceFracSurg}.}}
Since we assume that $S^3_{-(2n-1)/2}(C)$
is an $L$-space,  Proposition~\ref{prop:LargerLSpaces} guarantees that
$S^3_{-n}(C)$ is, too.
Our goal is to show that 
if $S^3_{-(2n-1)/2}(C)$ is an $L$-space,
then the inequalities from Lemma~\ref{lemma:Positivity} are all equalities:
\begin{equation}
\label{eq:LSpaceIneq}
d(S^3_{-(2n-1)/2}(C),v_i)-d(S^3_{-(2n-1)/2}(O),v_i)
= d(S^3_{-n}(C),w_i)-d(S^3_{-n}(O),w_i).
\end{equation}
for $i=0,...,2n-2$. 

By inspecting Equations~\eqref{eq:EvenList} and~\eqref{eq:OddList},
observe that the list $\{v_i\}_{i=0}^{2n-2}$ contains
each $\SpinC$ structure over $S^3_{-n}(C)$ twice, except for
one. Specifically, in the case where $n=2k$, the $\SpinC$ structure
which appears only once in the list is $v_0$. 
Combining  Lemma~\ref{lemma:Positivity} with the hypothesis
that $d(S^3_{-(2n-1)/2}(C),w_0)= d(S^3_{-(2n-1)/2}(O),w_0)$, we conclude that
\begin{equation}
\label{eq:ExtremeVanishEven}
\begin{array}{ll}
0=d(S^3_{-n}(C),v_0)-d(S^3_{-n}(O),v_0) & {\text{when $n=2k$.}}
\end{array}
\end{equation}
Similarly, when $n=2k+1$, the $\SpinC$ structure appearing only once on the
list is $v_{2k}$. It is easy to see that $v_{2k}={\overline{v_0}}$,
and hence combining the conjugation symmetry (Equation~\eqref{eq:ConjugationSymmetry}) with Lemma~\ref{lemma:Positivity}, we get that
\begin{equation}
\label{eq:ExtremeVanishOdd}
\begin{array}{ll}
0=d(S^3_{-n}(C),v_{2k})-d(S^3_{-n}(O),v_{2k}) &{\text{when $n=2k+1$.}}
\end{array}
\end{equation}

Note that Lemma~\ref{lemma:Positivity}
gives $2n-1$ inequalities which we must prove are all equalities.
Adding up these inequalities, we get
\begin{eqnarray}
\lefteqn{\sum_{i=0}^{2n-2}\left(d(S^3_{-(2n-1)/2}(C),w_i)-d(S^3_{-(2n-1)/2}(O),w_i)\right)} \nonumber \\
&\geq&
\sum_{i=0}^{2n-2}\left(d(S^3_{-n}(C),v_i)-d(S^3_{-n}(O),v_i)\right) \\
\label{ineq:CompareCorrTerms}
&=&
2 \sum_{\spinct\in\SpinC(S^3_{-n}(K))}
\left(d(S^3_{-n}(C),\spinct)-d(S^3_{-n}(O),\spinct)\right),
\nonumber
\end{eqnarray}
where here the last equation uses the fact that
every $\spinct\in\SpinC(S^3_{-n}(K))$ is represented twice
amongst the $\{v_i\}_{i=0}^{2n-2}$, except for one 
$v_j=v_0$ or $v_{2k}$, depending on the parity of $n$)
for which
we already know that 
$d(S^3_{-n}(C),v_j)-d(S^3_{-n}(O),v_j)=0$
(Equations~\eqref{eq:ExtremeVanishEven} and~\eqref{eq:ExtremeVanishOdd}).

Proposition~\ref{prop:OurEuler} now forces the inequality to be
equality; and this in turn implies that each of the $2n-1$ individual
inequalities in Inequality~\eqref{eq:dIneq} are also equalities.
\qed
\vskip.2cm

\vskip.2cm
\noindent{\bf Proof of Theorem~\ref{thm:LSpaceSymmetry}.}
By inspecting Equations~\eqref{eq:EvenList} and~\eqref{eq:OddList}, we
see that for $i=0,...,2n-2$ we have that
$c_1(w_i)/2=\phi(i)$. Moreover, another glance at those definitions
reveals that when $n=2k$, $v_i$ and $v_{2k-i}$ are conjugate $\SpinC$
structures for $i=1,...,k-1$.  Similarly, when $n=2k+1$, then $v_i$ and
$v_{2k-i}$ are conjugate for $i=0,...k-1$. It follows that
$$d(S^3_{-n}(C),v_i)-d(S^3_{-n}(O),v_i)=d(S^3_{-n}(C),v_{2k-i})-d(S^3_{-n}(O),v_{2k-i}).$$
Thus, the theorem follows from Theorem~\ref{thm:LSpaceFracSurg}.
\qed

\subsection{Euler characteristics}
\label{subsec:Euler}

We now return to the proof of Proposition~\ref{prop:OurEuler} stated
above. Indeed, we prove a more general statement. To give the
statement, we introduce some notation. If $Y$ is a rational homology
three-sphere and $k$ is some constant, then we define $$\HFp_{\leq
k}(Y)=\bigoplus_{\{d\in\Q\big| d\leq k\}}
\bigoplus_{\spinc\in\SpinC(Y)}\HFp_d(Y,\spinc).$$
In the case where $H_1(Y;\Z)\cong \Z$, we let $$\HFp_{\leq
k}(Y)=\bigoplus_{\{\spinc\in\SpinC(Y) \big|c_1(\spinc)\neq
0\}}\HFp(Y,\spinc)\oplus \bigoplus_{\{d\in\Q\big| d\leq k\}}
\HFp_d(Y,\spinc_0),$$ where $c_1(\spinc_0)=0$. 

In this latter case, the structure of $\HFinf$
(c.f. Equation~\eqref{eq:HFinfBOne}) ensures that for all sufficiently
large integers $N$, $\chi(\HF_{\leq 2N+\OneHalf}(Y))$ is independent
of $N$. We denote this constant by $\chiTrunc(\HFp(Y))$.

\begin{prop}
\label{prop:Euler}
For fixed relatively prime integers $p$ and $q$
with $p> 0$, there
is a constant $k(p,q)$ with the property that for any knot $C\subset
S^3$, $$\sum_{i\in\Zmod{p}} \left( \chi
\left(\HFpRed(S^3_{p/q}(C),i)\right)-
\frac{d(S^3_{p/q}(C),i)}{2}\right)
- q \cm \chiTrunc(\HFp(S^3_0(C)))=k(p,q)$$
\end{prop}

We break the proof into two steps.

\begin{lemma}
\label{lemma:EulerTruncRed}
Let  $p$ and $q$ be 
relatively prime integers with $p\geq 0$. There is a constant
$k_1(p,q)$ with the property that for any knot $C\subset S^3$
and any sufficiently large integer $N$,
\begin{eqnarray}
\label{eq:EulerLemma}
\lefteqn{\chi(\HFp_{\leq 2N}(S^3_{p/q}(C)))-N\cm p}&& \\
&=& \sum_{\spinc\in\SpinC(S^3_{p/q}(C))}
\left(
\chi(\HFpRed(S^3_{p/q}(C))-\frac{d(S^3_{p/q}(C),\spinc)}{2}\right) \nonumber + k_1(p,q). \nonumber
\end{eqnarray}
\end{lemma}

\begin{proof}
Let $Y=S^3_{p/q}(C)$.
For sufficiently large $N$, $\HFpRed(Y)$ is contained in $\HFp_{\leq 2N}(Y)$.
Over $\Field$, we have a splitting
$$
\HFp_{\leq 2N}(Y)\cong \HFpRed(Y)\oplus (\Image\pi \cap \HFp_{\leq 2N}(Y)).
$$
But it follows readily from the structure of $\HFinf(Y)$
(c.f. Equation~\eqref{eq:HFinf}) that
\begin{eqnarray*}
\chi(\Image\pi\cap \HFp_{\leq 2N}(Y)) &=& \sum_{\spinc\in\SpinC(Y)}
\#\{{[d(Y,\spinc),2N]}\cap (d+2\Z)\subset \Q\} \\
&=& \sum_{\{\spinc\in\SpinC(Y)\}} \left(N + 1 -\lceil \frac{d(Y,\spinc)}{2}\rceil\right),
\end{eqnarray*}
where here $\lceil x\rceil$ denotes the smallest integer greater than or
equal to $x$.  Thus, since the number of $\SpinC$ structures on $Y$ is
$p$, we get that
\begin{eqnarray*}
\lefteqn{\chi(\HFp_{\leq 2N}(Y))-N\cm p} \\
 &=&
 \sum_{\spinc\in\SpinC(Y)} \left(\HFpRed(Y,\spinc)-\lceil\frac{d(Y,\spinc)}{2}\rceil+1\right) \\
&=&
\sum_{\spinc\in\SpinC(Y)} \left(\HFpRed(Y,\spinc)-\frac{d(Y,\spinc)}{2}\right) 
+ \sum_{\spinc\in\SpinC(Y)}
\left(\frac{d(Y,\spinc)}{2}-\lceil\frac{d(Y,\spinc)}{2}\rceil+1\right).
\end{eqnarray*}
but this last sum is easily seen to depend on $d(Y,\spinc)$ only modulo $2\Z$; i.e.
it depends only on the linking form of $Y=S^3_{p/q}(C)$, which in turn is independent
of the knot $C$ (depending only on $p$ and $q$).
\end{proof}

\begin{lemma}
\label{lemma:EulerTrunc}
Given a pair of non-negative, relatively prime integers
$p$ and $q$, there is a constant $k_2(p,q)$ with the property
that for any knot $C\subset S^3$, there is a natural number $N_0$ with
the property that for all $N\geq N_0$, $$\chi(\HFp_{\leq
2N}(S^3_{p/q}(C))) -N\cm p = q \cm \chiTrunc(\HFp(S^3_0(C))) +
k_2(p,q).$$
\end{lemma}

\begin{proof}
We will use induction on $p+q$, together with the fact that 
for any pair 
$(p,q)$ of relatively
prime, non-negetive integers with $p+q>1$, 
there are two pairs of non-negative,
relatively integers $(p_0,q_0)$ and $(p_2,q_2)$ with the properties
that 
\begin{eqnarray}
p_0 \cm q - p \cm q_0 &=& -1 \label{eq:RelPrime} \\
(p,q)&=& (p_0,q_0)+(p_2,q_2) \label{eq:SumProperty}
\end{eqnarray}

In the basic case of the lemma 
where $p+q=1$, we consider cases where $(p,q)=(1,0)$
or $(0,1)$. In both cases, the lemma is clear.

For arbitrary relatively prime, non-negative integers $(p,q)$, find
non-negative integers $(p_0,q_0)$ and $(p_2,q_2)$ satisfying
Equations~\eqref{eq:RelPrime} and~\eqref{eq:SumProperty}.  Let
$Y_0=S^3_{p_0/q_0}(C)$, $Y_1=S^3_{p/q}(C)$, and
$Y_2=S^3_{p_2/q_2}(C)$. Those equations guarantee that we have a long
exact sequence $$\begin{CD} ...@>>> \HFp(Y_0)
@>{f_0}>>\HFp(Y_1) @>{f_2}>>\HFp(Y_2) @>{f_3}>> ...,
\end{CD}$$
and also that the lemma is known for $(p_0,q_0)$ and $(p_2,q_2)$.
We assume first that that $p_0\neq 0$. 

When $N$ is sufficiently large, the restriction $g_0$ of $f_0$ to
$\HFp_{\leq 2N}(Y_0)$ is contained in $\HFp_{\leq
2N+\frac{1}{4}}(Y_1)$, the restriction of $f_2$ to $\HFp_{\leq
2N+\frac{1}{4}}(Y_1)$ is contained in $\HFp_{\leq 2N+\OneHalf}(Y_2)$,
and the finally, the restriction of $f_3$ to $\HFp_{\leq
2N+\OneHalf}(Y_2)$ is contained in $\HFp_{\leq 2N}(Y_0)$. This follows
at once from the grading shift formula,
Equation~\eqref{eq:GradingShift}: we have that
$\chi(W_i)=1$ and $\sigma(W_i)=-1$ for $i=0,1$; while the cobordism
$W_2$ induces the trivial map on $\HFinf$ since $b_2^+(W_0)=1$.

Choosing $N$ as above, consider the diagram
$$
\begin{CD}
&& 0 &&  0 && 0  \\
&& @VVV @VVV @VVV \\
... @>{g_2}>> \HFp_{\leq 2N}(Y_0) @>{g_0}>> \HFp_{\leq 2N+\frac{1}{4}}(Y_1) @>{g_1}>> \HFp_{\leq 2N+\OneHalf}(Y_2) @>{g_2}>> ... \\
&& @VVV @VVV @VVV \\
... @>{f_2}>> \HFp(Y_0) @>{f_0}>>\HFp(Y_1) @>{f_1}>> \HFp(Y_2) @>{f_2}>> ... \\
&& @VVV @VVV @VVV \\
... @>{h_2}>> \HFp_{>2N}(Y_0)@>{h_0}>>\HFp_{>2N+\frac{1}{4}}(Y_1) @>{h_1}>>\HFp_{>2N+\OneHalf}(Y_2)  @>{h_2}>> ... \\
&& @VVV @VVV @VVV \\
&& 0 &&  0 && 0,  \\
\end{CD}
$$ where the columns are exact. Note that the rows are not necessarily
exact (except for the middle one), however all three rows can be
thought of as chain complexes. We denote these three rows by ${\mathcal R}_1$,
${\mathcal R}_2$, and ${\mathcal R}_3$. Since ${\mathcal R}_2$ is exact, it follows that $H_*({\mathcal R}_1)\cong H_*({\mathcal R}_3)$.

We claim that $H_*({\mathcal R}_3)$ is independent of $C$ and $N$
(provided the latter is sufficiently large). To see this, observe that
we have a diagram $$
\begin{CD}
... @>{h_2^{\infty}}>> \HFinf_{>2N}(Y_0)@>{h_0^\infty}>>\HFinf_{>2N+\frac{1}{4}}(Y_1) @>{h_1^\infty}>>\HFinf_{>2N+\OneHalf}(Y_2)  @>{h_0^{\infty}}>> ... \\
&& @V{\cong}VV @V{\cong}VV @V{\cong}VV \\
... @>{h_2}>> \HFp_{>2N}(Y_0)@>{h_0}>>\HFp_{>2N+\frac{1}{4}}(Y_1) @>{h_1}>>\HFp_{>2N+\OneHalf}(Y_2)  @>{h_0}>> ..., \\
\end{CD}
$$
where here $h_0$ is the sum over all $\spinc\in\SpinC(W_0)$ of the projections of the
induced maps on $\HFinf$; e.g. letting 
$$\Pi_{>2N+\OneHalf}\colon \HFinf(Y_1) \longrightarrow \HFinf_{>2N+\OneHalf}(Y_1)$$
denote the projection, we let $h_0^{\infty}$ be the restriction to $\HF_{>2N}(Y_0)$ of 
$$\sum_{\spinc\in\SpinC(W_0)} \Pi_{>2N+\OneHalf}\circ F^{\infty}_{W_0,\spinc}.$$
The map $h_i^{\infty}$ are defined similarly. Note that
$h_2^{\infty}=0$, since the map induced by $W_2$ has $b_2^+(W_2)=1$
(in general, this is proved in
(c.f. Lemma~\ref{HolDiskFour:lemma:BTwoPlusLemma}
of~\cite{HolDiskFour}, though in the present case it follows
quickly from the dimension formula: the cobordism 
$W_2$ reverses the $\Zmod{2}$ grading, but $\HFinf$ is supported entirely
in even gradings).

So far we have established that for all sufficiently large $N$,
$$\chi (H_*({\mathcal R}_1)) = \chi (\HFp_{\leq 2N}(Y_0)) - \chi(\HFp_{\leq 2N+\OneQuarter}(Y_1))
+\chi(\HFp_{\leq 2N+\OneHalf}(Y_2))$$
is independent of $C$ and $N$ (provided that the latter is sufficiently large);
but it is also clear that for sufficiently large $N$, 
\begin{eqnarray*}
\chi(\HFp_{\leq 2N+\OneQuarter}(Y_1))=\chi (\HFp_{\leq 2N}(Y_1)) + c_3 &{\text{and}}&
\chi(\HFp_{\leq 2N+\OneHalf}(Y_2))=\chi (\HFp_{\leq 2N}(Y_2)) + c_4,
\end{eqnarray*}
with constants $c_3$ and $c_4$ again depending only on $(p_0,q_0)$ and $(p_2,q_2)$ respectively.
Combining all the constants, we establish the inductive step, at least in
the case where $p_0$ is non-zero.

In the case where  $p_0=0$, the above argument works with
slight modification. In this case, we cyclically order the three
three-manifolds $(Y_0, Y_1, Y_2)$ so that $Y_1$ has $b_1(Y_1)=1$. In
this case, the dimension shifts works slightly differently:
$\sigma(W_0)=\sigma(W_1)=0$ and hence, we compare $\HFp_{\leq
2N}(Y_0)$, $\HFp_{\leq 2N+\OneHalf}(Y_1)$, and $\HFp_{\leq 2N+1}(Y_2)$
(to see that the maps induced by $W_2$ carry $\HFp_{\leq 2N+1}(Y_2)$
for sufficiently large $N$, note that Equation~\eqref{eq:GradingShift}
ensures that the the shift in grading is at least zero, and also that
$\HFp_{\leq 2N+1}(Y_2)=\HFp_{\leq 2N}(Y_2)$ for sufficiently large $N$:
$\HFinf$ is supported only in even degrees). With these minor
modifications, the previous argument establishes the inductive step in
the remaining case, as well.
\end{proof}

\vskip.2cm
\noindent{\bf{Proof of Proposition~\ref{prop:Euler}.}}
When $p$ and $q$ are non-negative,
this is a combination of Lemmas~\ref{lemma:EulerTruncRed} and
\ref{lemma:EulerTrunc}. This could be proved either 
by running the induction from Lemma~\ref{lemma:EulerTrunc} to show
that Lemma~\ref{lemma:EulerTrunc} still holds in the case where $p>0$
and $q\leq 0$; or alternatively, deducing this from the previous case,
using the facts that $S^3_{p/q}(C)=-S^3_{-p/q}(r(C))$ (where $r$
denotes reflection), $\chi(\HFpRed(-Y))=-\chi(\HFpRed(Y))$, and the
fact that $\chiTrunc(\HFp(S^3_0(C)))=\chiTrunc(\HFp(S^3_0(r(C))))$,
all of which are established in~\cite{HolDiskTwo}.
\qed
\vskip.2cm

\vskip.2cm
\noindent{\bf{Proof of Proposition~\ref{prop:OurEuler}.}}
Note that if $S^3_{-(2n-1)/2}(C)$ is an $L$-space,
then according to Proposition~\ref{prop:LargerLSpaces}, so
is $S^3_{-(2n-1)/2}(C)$, and in particular
$\chi(\HFpRed(S^3_{-(2n-1)/2}(C),i))=0$
and $\chi(\HFpRed(S^3_{-n}(C),i))=0$. The rest is a direct 
application of Proposition~\ref{prop:Euler}.
\qed
\vskip.2cm

\section{Proof of Theorem~\ref{thm:AlternatingCriterion}.}
\label{sec:Proof}

\vskip.2cm
\noindent{\bf{Proof of Theorem~\ref{thm:AlternatingCriterion}.}}
If $K$ has unknotting number one, by Montesinos' trick
(Lemma~\ref{lemma:Montesinos}), we know that $\Sigma(K)=\pm
S^3_{D/2}(C)$. After reflecting $K$ if necessary, we can achieve that
$\Sigma(K)=S^3_{-D/2}(C)$.  According to
Proposition~\ref{prop:AltKnots}, $\Sigma(K)$ is an $L$-space, and we
have an isomorphism $$\phi\colon\Coker
(q)\stackrel{\cong}{\longrightarrow} H^2(\Sigma(K);\Z) $$ with the
property that $$\MinVecLen_Q(\xi)=d(\Sigma(K),\phi(\xi)).$$
(Note that $\Sigma(r(K))=-\Sigma(K)$, and hence the reflection
has the effect of reversing the signs of the correction terms, c.f.
Equation~\eqref{eq:ReverseOrientation};
this is responsible for the sign $\epsilon$ appearing in the statement
of the theorem.)

Now the expression of $S^3_{-D/2}\cong \Sigma(K)$ gives us an
identification $\Zmod{D} \cong \Coker(q)$. Thus, according to
Theorem~\ref{thm:IntForms}, we get both
\begin{eqnarray*}
\gamma_{D}(i)&{\equiv}& \MinVecLen_Q(\phi(i))\pmod{2\Z} \\
\gamma_{D}(i)&\leq& \MinVecLen_Q(\phi(i)),
\end{eqnarray*}
where as usual $\gamma_{D}(i)=\MinVecLen_{R_D}(i)$.

We claim that $\MinVecLen_{R_D}(i)=d(S^3_{-D/2}(O),i)$. This can be seen,
for example, by taking a two-bridge knot whose branched double-cover is
$S^3_{-D/2}(0)$, and applying Proposition~\ref{prop:AltKnots}.

Note that for $D=2n-1$,
$$\MinVecLen_{R_D}(0)=\left\{\begin{array}{ll}
0 & {\text{if $n$ is odd}} \\
\OneHalf & {\text{if $n$ is even.}}
\end{array}
\right.$$
Thus, when $|\MinVecLen_{Q_K}(0)|\leq \OneHalf$, it
follows that $\MinVecLen_{R_D}(0)=\MinVecLen_{Q_K}(0)$, verifying the
hypothesis of Theorem~\ref{thm:LSpaceSymmetry}.  The symmetry
$T_{\phi}(i)=T_{\phi}(2k-i)$ is now a direct consequence of 
Theorem~\ref{thm:LSpaceSymmetry}.  
\qed
\vskip.2cm

We note that the obstruction given by
Theorem~\ref{thm:AlternatingCriterion} does not depend on the choice
of the alternating projection of $K$, since, according to 
Proposition~\ref{prop:AltKnots},  $M_Q$ is a topological
invariant of the oriented branched double-cover. Note also
that if we reflect $K$
(or, equivalently, use the black instead of the white graphs), this
has the effect of reversing the orientation of $\Sigma(K)$, and hence
the corresponding correction terms (as given by $M_Q$ for the Goeritz
form, c.f. Proposition~\ref{prop:AltKnots}) all get multiplied by
$-1$. This is compensated by our freedom in choosing $\epsilon=\pm 1$.
We return to this point in Section~\ref{sec:Refine}.

\section{Calculations for alternating knots}
\label{sec:Calcs}

We explain now how to apply Theorem~\ref{thm:AlternatingCriterion} in detail.

Given an alternating knot of determinant $D$, we start by writing down
its $m\times m$ Goeritz matrix $G$. (In practice, it can be useful to
reflect the knot if necessary to minimize the number of white
regions.) Suppose that $H_1(\Sigma(K);\Z)$ is cyclic or, equivalently,
that $\Coker(q)\cong\Zmod{D}$.

Next, we find the function $\MinVecLen_Q\colon \Coker(q)\longrightarrow
\Q$.  Two vectors $v_1, v_2\in \Z^m$ correspond to equivalent vectors
in $\Coker(q)$ when $G^{-1}(v_1-v_2)\in\Z^m$. The induced quadratic
form $Q^*$ in this basis is represented by $(v,w) \mapsto v^t\cm
G^{-1} \cm w$. We claim that characteristic vectors in $V^*\cong \Z^m$
which achieve maximal length are contained in the finite set
$$\{v=(v_1,...,v_m)\in\Z^m\big|~{\text{for $i=1,...,m$,}}~ |v_i|\leq |G_{i,i}|~\text{and}~ v_i\equiv G_{i,i}\pmod{2}\};$$ if a
vector lies outside that set, it is straightforward to find another
equivalent vector $v'$ with larger length. 

Thus, by performing a finite set of calculations, we end up with a
list of $D$ vectors in $\Z^m$ which maximize their length in their
equivalence class. We order these vectors $\{x{(i)}\}_{i=0}^{D-1}$, so
that $x{(i)}$ represents $i\in\Zmod{D}$ under some isomorphism of
$\Coker(q)\cong\Zmod{D}$. Next, form the vector $\{a_i\}_{i=0}^{D-1}$,
where $$4a_i=x{(i)}^{t}\cm G^{-1}\cm x{(i)}+m.$$ Note that this
ordering of the vectors in $A$ is not canonical: it is canonical only
up to reordering given by multiplication by the units in
$\Zmod{D}$. According to Proposition~\ref{prop:AltKnots}, the vector
$A$ contains the correction terms for $\Sigma(K)$.

Next, we calculate $B_i=\gamma_D(i)$. This is straightforward: let
$4B_i=\kappa(i)^{t}\cm R_{D}^{-1}\cm \kappa(i)+2$, where the
$\kappa(i)$ are given in order in Equations~\eqref{eq:EvenList} and
\eqref{eq:OddList}. 

Now for each automorphism $\phi$ of $\Zmod{D}$ and sign $\epsilon=\pm
1$, we form the vector $C_i=-B_i-\epsilon A_{\phi(i)}$. We call such a
vector $C_i$ a {\em matching} for the knot $K$. The set of matchings
for an alternating knot $K$ is a knot invariant. We call a matching
{\em even} if it consists of even integers. We call a matching {\em
positive} it consists of non-negative rational numbers. Finally,
writing $D=4k\pm 1$, we call a matching {\em symmetric} if it
satisfies the symmetry $C_{i}=C_{2k-i}$ for $1<i<k$ and also for $i=0$
when $D=4k+1$. Note that matchings always satisfy the symmetry
$C_{i}=C_{D-i}$, and hence the matching is determined by
$\{C_i\}_{i=0}^{n-1}$, where $D=2n-1$.

In this language, Theorem~\ref{thm:AlternatingCriterion} says that
if $K$ is an alternating knot with unknotting number equal to one, then
there is at least one matching $C$ which is positive and even. 
Moreover, if $|A_0|\leq \OneHalf$, then there is also an
even, positive, and symmetric matching. (Note that the
condition that $|A_0|\leq \OneHalf$ is equivalent to the existence
of some matching for which $C_0=0$.)
Note also that all the knots we consider in this section
satisfy the condition that $|A_0|\leq \OneHalf$.

We now turn to the applications of Theorem~\ref{thm:AlternatingCriterion}
to knots with small crossing numbers. 

\vskip.2cm
\noindent{\bf{Proof of Corollary~\ref{cor:VerySmallKnots}.}}

\begin{figure}
\mbox{\vbox{\epsfbox{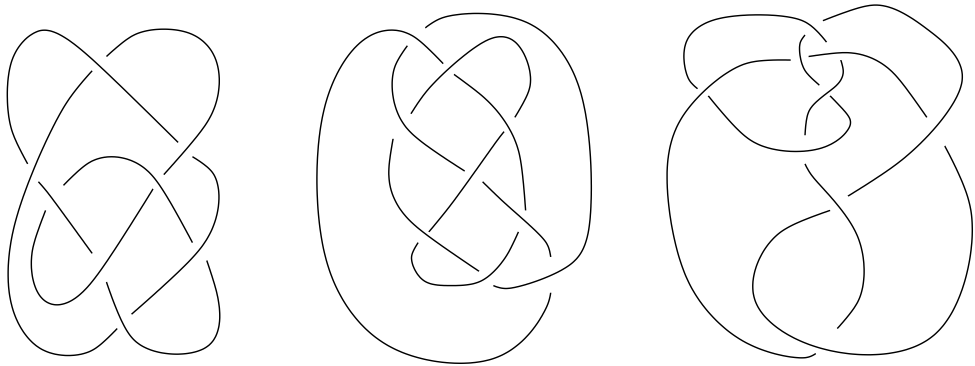}}}
\caption{\label{fig:SmallKnots}
{\bf The knots $8_{10}$, $9_{29}$, and $9_{32}$.}}
\end{figure}

We start with the knot $8_{10}$ pictured in Figure~\ref{fig:SmallKnots}.
Its Goeritz form is represented by the matrix
$$ 
G=\left(\begin{array}{rrr}
-4,     &       1,      &       1       \\
1,      &       -2,     &       1       \\
1,      &       1,      &       -5
\end{array}\right),$$
whose determinant is $-27$.

%       \begin{figure}
%       \mbox{\vbox{\epsfbox{8s10.eps}}}
%       \caption{\label{fig:8s10}
%       {\bf The knot $8_{10}$.}}
%       \end{figure}

Following the above procedure, we find the maximal squares of
lengths of the
vectors in $\Z^3$ in each equivalence class or, more precisely, divide
these numbers by four, add $3/4$, and then order according to the
group structure of $\Zmod{27}$.  This gives us the ordered list of
numbers: $$ A = \left(
\begin{array}{rrrrrrrrr}
- \frac{1}{2}  &\frac{25}{54} &
   - \frac{35}{54}  &\frac{1}{6} &- \frac{59}{54}  &
   - \frac{23}{54}  &\frac{1}{6} &\frac{37}{54} &
   - \frac{47}{54} \\ &&&&&&&& \\- \frac{1}{2}  &
   - \frac{11}{54}  &\frac{1}{54} &\frac{1}{6} &\frac{13}{54} &
   \frac{13}{54} &\frac{1}{6} &\frac{1}{54} &- \frac{11}{54}   \\ &&&&&&&&  \\
   - \frac{1}{2}  & 
  - \frac{47}{54}  &\frac{37}{54} &
   \frac{1}{6} &- \frac{23}{54}  &- \frac{59}{54}  &
   \frac{1}{6} &- \frac{35}{54}  &\frac{25}{54}
\end{array}
\right).$$

We have written these in an order compatible with the isomorphism 
of $\Coker(q)\cong \Zmod{27}$, where the first
term corresponds to the spin structure (and the $i^{th}$ term corresponds
to the cohomology class $(0,0,i)$.

We compare this ordered list with the list
$$
B = 
 \left(
\begin{array}{rrrrrrrrr}
\frac{1}{2} & \frac{23}{54} & \frac{11}{54} & - \frac{1}{6}   & 
   - \frac{37}{54}   & - \frac{73}{54}   & 
   - \frac{13}{6}   & - \frac{169}{54}   & 
   - \frac{121}{54} \\ &&&&&&&& \\
    - \frac{3}{2}   & 
   - \frac{49}{54}   & - \frac{25}{54}   & 
   - \frac{1}{6}   & - \frac{1}{54}   & 
   - \frac{1}{54}   & - \frac{1}{6}   & 
   - \frac{25}{54}   & - \frac{49}{54}     \\ &&&&&&&& \\
   - \frac{3}{2}   & - \frac{121}{54}   & 
   - \frac{169}{54}   & - \frac{13}{6}   & 
   - \frac{73}{54}   & - \frac{37}{54}   & 
   - \frac{1}{6}   & \frac{11}{54} & \frac{23}{54},
\end{array}
\right). 
$$
which are $\gamma_{27}(i)$ for $i=0,...,26$.

By comparing $A_0$ and $B_0$ we see that an even matching (for some
automorphism $\phi$ and $\epsilon=\pm 1$) can exist only when
$\epsilon=1$.  By inspection, we see that there are only two possible
automorphisms $\phi$ of $\Zmod{27}$ (multiplication by $\pm 5$) for
which $C_i=-B_i-A_{\phi(i)}$ is a non-negative, even integer for
$i=0,...,26$.  For both of these, $C_i=-B_i-A_{\phi(i)}$
($T_{1,\phi}$) is given by the list $$ 0, 0, 0, 0, 0, 2, 2, 4, 2, 2,
2, 0, 0, 0, 0, 0, 0, 2, 2, 2, 4, 2, 2, 0, 0, 0, 0.  $$ The matching
evidently fails the symmetry: $T_{1,\phi}(4)\neq T_{1,\phi}(10)$.  So
there are no even, positive, symmetric matchings.  Thus, by
Theorem~\ref{thm:AlternatingCriterion}, $8_{10}$
cannot have unknotting number equal to one.

We abbreviate this data somewhat.
If $C$ is a matching, we will list only the first $n$ terms 
(where $D=2n-1$), dropping all the initial and final zero terms, and
indicating the ${k}^{th}$ term in bold face.  
For example, we indicate the even, positive matching
for $8_{10}$ by
$$2, 2, {\bf 4}, 2, 2, 2.$$ The fact that this matching
is asymmetric is now obvious from this notation.

For $9_{29}$, there is only one possible choice of $\epsilon$ for
which we can find a $\phi$ satisfying
Condition~\eqref{eq:Integrality}; and for that choice, there are four
possible $\phi$ satisfying Conditions~\eqref{eq:Integrality} and
\eqref{eq:Positivity}, all of which give the same 
matching: $$ 2, 2, 2, 2,
4, 4, 6, {\bf 6}, 6, 4, 4, 4, 2, 2, 2, 2.  $$ For $9_{32}$, again
there is only one even, positive 
matching: $$ 2, 2, 2, 2,
4, 4, 6, 6, {\bf 8}, 6, 6, 4, 4, 4, 2, 2, 2, 2. $$ Since both of 
the above
matchings are not symmetric, Theorem~\ref{thm:AlternatingCriterion}
shows that the knots do not have unknotting number equal to one.  

On the other hand, it is clear from their pictures that the three
knots considered here can be unknotted in two steps. \qed

\subsection{Classification of $\leq 9$-crossing knots with $u=1$}

Note that this classification already follows from
Corollary~\ref{cor:VerySmallKnots}, together with previously known
results as explained in the introduction. However, for completeness we
include a classification here using Murasugi's bound, the
cyclicity of $H_1(\Sigma(K);\Z)$,
and Theorem~\ref{thm:AlternatingCriterion}. 

First observe that the following knots fail Condition~\eqref{eq:Integrality} in Theorem~\ref{thm:AlternatingCriterion},
(i.e. they admit no even matchings, in the terminology from the beginning of this section):
\begin{equation}
\label{eq:QList}
7_4, 8_{8}, 8_{16}, 9_{15}, 9_{17}, 9_{31}.
\end{equation}
For example,
a Goeritz matrix for $8_{16}$ is given by
$$G=
\left(\begin{array}{rrr}
-4 & 1 & 1 \\
1 & -4 & 1 \\
1 & 1 & -3
\end{array}\right).
$$
The list of correction terms $A$ and also the vector $B$ respectively are listed as follows:
\[
\begin{tiny}
\left(\begin{array}{rrrrrrr}
 - \frac{1}{2}  &- \frac{43}{70}  &
   - \frac{67}{70}  &\frac{33}{70}&
   - \frac{23}{70}  &\frac{9}{14}&- \frac{43}{70} \\ \\
   - \frac{1}{10}  &\frac{13}{70}&\frac{17}{70}&\frac{1}{14}&
   - \frac{23}{70}  &- \frac{67}{70}  &
   \frac{13}{70} \\ \\ - \frac{9}{10}  &- \frac{3}{14}  &
   \frac{17}{70}&\frac{33}{70}&\frac{33}{70}&\frac{17}{70}&
   - \frac{3}{14} \\  \\- \frac{9}{10}  &\frac{13}{70}&
   - \frac{67}{70}  &- \frac{23}{70}  &\frac{1}{14}&
   \frac{17}{70}&\frac{13}{70} \\ \\- \frac{1}{10}  &
   - \frac{43}{70}  &\frac{9}{14}&- \frac{23}{70}  &
   \frac{33}{70}&- \frac{67}{70}  &
   - \frac{43}{70}  
\end{array}\right)
\left(
\begin{array}{rrrrrrr}
 \frac{1}{2}&- \frac{289}{70}  &
   - \frac{1}{70}  &- \frac{221}{70}  &
   \frac{31}{70}&- \frac{45}{14}  &- \frac{9}{70}  \\ \\
   - \frac{23}{10}  &\frac{19}{70}&
   - \frac{169}{70}  &- \frac{5}{14}  &
   - \frac{109}{70}  &- \frac{1}{70}  &
   - \frac{121}{70}  \\ \\ - \frac{7}{10}  &
   - \frac{13}{14}  &- \frac{29}{70}  &
   - \frac{81}{70}  &- \frac{81}{70}  &
   - \frac{29}{70}  &- \frac{13}{14}  \\ \\
   - \frac{7}{10}  &- \frac{121}{70}  &
   - \frac{1}{70}  &- \frac{109}{70}  &
   - \frac{5}{14}  &- \frac{169}{70}  &
   \frac{19}{70} \\ \\- \frac{23}{10}  &- \frac{9}{70}  &
   - \frac{45}{14}  &\frac{31}{70}&
   - \frac{221}{70}  &- \frac{1}{70}  &
   - \frac{289}{70}  
\end{array}
 \right).
\end{tiny}
\]
By looking at the $A_{7i}$ and $B_{7j}$ terms, we see that there is no even matching.
Similar calculations rule out the other knots in the above list.
%       This is a straightforward calculation. For example, 
%       the correction terms for $7_4$ are given by the list:
%       $$  - \frac{1}{2}  ,\frac{17}{30},- \frac{7}{30}  ,
%          \frac{11}{10},\frac{17}{30},\frac{1}{6},- \frac{1}{10}  ,
%          - \frac{7}{30}  ,- \frac{7}{30}  ,
%          - \frac{1}{10}  ,\frac{1}{6},\frac{17}{30},\frac{11}{10},
%          - \frac{7}{30}  ,\frac{17}{30}$$
%       while $\gamma_{15}$ is given by the list
%       $$   \frac{1}{2},\frac{11}{30},- \frac{1}{30}  ,
%          - \frac{7}{10}  ,- \frac{49}{30}  ,
%          - \frac{5}{6}  ,- \frac{3}{10}  ,
%          - \frac{1}{30}  ,- \frac{1}{30}  ,
%          - \frac{3}{10}  ,- \frac{5}{6}  ,
%          - \frac{49}{30}  ,- \frac{7}{10}  ,
%          - \frac{1}{30}  ,\frac{11}{30}.$$

The knot $9_5$ cannot have unknotting number equal to one since 
it has no even, positive matchings.
Indeed, there is only one even matching,
and it is given by
$$-2, 0, 0, 0, 2, {\bf 2}, 2.$$
In particular, this shows (following Theorem~\ref{thm:IntFormBounds})
that the branched double-cover $\Sigma(9_5)$ does not bound any four-manifold with intersection
form $R_{23}$.
(Of course, the branched double-covers of the 
knots in List~\eqref{eq:QList} also do not bound corresponding intersection forms, but 
for more elementary reasons.)

To rule out the remaining knots, we must use all three conditions in Theorem~\ref{thm:AlternatingCriterion}. Specifically, for each of the knots there are unique even, positive matchings,
but none of them are symmetric. The matchings are listed in the following
table:
\vskip.2cm
\begin{tabular}{llllllllllllll}
$8_{3}$:  &  &  &  &  & 2& {\bf 2}&  &  &  &  &  &  \\
$8_{4}$:  &  &  &  &  & 2& {\bf 2}& 2& 2&  &  &  &  \\
$8_{6}$:  &  &  &  &  & 2& {\bf 2}& 2& 2&  &  &  &  \\
$8_{12}$: &  &  & 2& 2& 2& {\bf 4}& 2& 2&  &  &  &  \\
$9_{8}$:  &  &  & 2& 2& 2& {\bf 4}& 4& 2& 2& 2&  &  \\
$9_{25}$: &2 & 2& 2& 4& 4& {\bf 6}& 4& 4& 4& 2& 2& 2
\end{tabular}
\vskip.2cm

\subsection{Corollary~\ref{cor:Alt10}}

Again, the knots appearing in the list for Corollary~\ref{cor:Alt10}
fail Theorem~\ref{thm:AlternatingCriterion} at several different
levels. For example, several of them cannot have unknotting number
equal to one, because they have no even matchings:
\begin{equation}
\label{eq:10altRat}
\begin{array}{lllllll} 10_{67} & 10_{86} & 10_{105} & 10_{106} &
10_{109} & 10_{116} & 10_{121}.
\end{array}
\end{equation}

The knot $10_{68}$ has no even, positive matchings:
We list its only even matching:
\[
\begin{array}{lc}
10_{68}: & 2, 2, 2, 4, 4, 4, 6, {\bf 6}, 6, 4, 4, 4, 2, 2, 2, 0, 0, 0,
 0, 0, -2.
\end{array}
\]
Again note that the branched double cover of the knot $10_{68}$ does
not bound a smooth four-manifold with intersection form of type $R_{D}$
(and also the same holds in the topological category for the knots
appearing in List~\eqref{eq:10altRat}).

For the remaining knots listed in Corollary~\ref{cor:Alt10}, we use
all three conditions. Below we list all the even, positive matchings
for the remaining knots.  Failure of the symmetry condition is now
evident. Note that the knots $10_{58}$ and $10_{77}$ appear twice on
this list, since they have two distinct matchings.
\vskip.2cm
\tabcolsep=0.8mm
\begin{small}
\begin{tabular}{ccccccccccccccccccccccccccccccccccccccc}
$10_{48}$:  &  &  &  &  &  &  &  &  &  &  & 2& 2& 2& 2&  4&  4&  4& {\bf 6}& 6& 4& 2& 2& 2& 2&  \\
$10_{51}$:  &  &  &  &  &  &  &  &  & 2& 2& 2& 4& 4& 4&  6&  6&  8& {\bf 8}& 6& 6& 6& 4& 4& 4& 2& 2& 2& 2& 2 \\
$10_{52}$:  &  &  &  &  &  &  &  &  & 2& 2& 2& 2& 4& 4&  4&  6&  6& {\bf 8}& 6& 6& 4& 4& 4& 2& 2& 2& 2& 0& 2& \\
$10_{54}$:  &  &  &  &  &  &  &  &  &  &  &  & 2& 2& 2&  4&  4&  6& {\bf 6}& 4& 4& 4& 2& 2& 2& 2& 2 \\
$10_{57}$:  &  &  &  &  &  & 2& 2& 2& 2& 2& 4& 4& 6& 6&  6&  8&  8& {\bf 10}& 8& 8& 8& 6&  6& 4&  4& 4& 2& 2& 2& 2& 2& \\
$10_{58}$:  &  &  &  &  &  &  &  &  & 2& 2& 2& 4& 4& 4&  6&  6&  6& {\bf 8}& 6& 6& 4& 4& 4& 2& 2& 2 \\
$10_{58}'$: &  &  &  &  &  &  & 2& 2& 2& 2& 2& 4& 4& 4&  6&  6&  6& {\bf 8}& 8& 6& 6& 4& 4& 4& 2& 2&  2 \\
$10_{64}$:  &  &  &  &  &  &  &  &  &  &  &  & 2& 2& 4&  4&  4&  6& {\bf 6}& 6& 4& 4& 2& 2& 2& 2& 2& \\
$10_{70}$:  &  &  &  &  &  &  &  & 2& 2& 2& 2& 4& 4& 4&  6&  6&  8& {\bf 8}& 8& 6& 6& 4& 4& 4& 2& 2& 2& 2& 0& 2& \\
$10_{77}$:  &  &  &  &  &  &  &  &  &  & 2& 2& 2& 4& 4&  6&  6&  6& {\bf 8}& 6& 6& 4& 4& 4& 2& 2& 2&  2& 2 \\ 
$10_{77}'$: &  &  &  &  &  &  &  &  &  & 2& 2& 2& 4& 4&  4&  6&  6& {\bf 8}& 6& 6& 6& 4& 4& 4& 2& 2&  2& 0 \\
$10_{79}$:  &  &  &  &  &  &  &  & 2& 0& 2& 2& 2& 4& 4&  4&  6&  6& {\bf 8}& 6& 6& 4& 4& 4& 2& 2& 2 \\
$10_{81}$:  &  &  &  & 2& 2& 2& 2& 4& 4& 4& 4& 6& 6& 6&  8&  8& 10& {\bf 10}& 10& 8& 8& 6& 6& 4& 4& 4& 4& 2& 2& 2& 2 \\
$10_{83}$:  &  &  &  &  & 2& 2& 2& 2& 2& 4& 4& 4& 6& 6&  8&  8&  8& {\bf 10}& 10& 8&  8& 6& 6& 6& 4& 4& 4& 2& 2& 2& 2& 2  \\
$10_{87}$:  &  &  &  & 2& 2& 2& 2& 2& 4& 4& 4& 6& 6& 6&  8&  8& 10& {\bf 10}& 8& 8& 8&   6& 6& 4& 4& 4& 2& 2& 2& 2& 2 \\
$10_{90}$:  &  &  &  & 2& 0& 2& 2& 2& 2& 4& 4& 4& 6& 6&  6&  8&  8& {\bf 10}& 8& 8& 6& 6& 6&   4& 4& 4& 2& 2& 2& 2 \\
$10_{93}$:  &  &  &  &  &  &  &  & 2& 2& 2& 2& 4& 4& 6&  6&  6&  8& {\bf 8}& 8& 6& 6& 4& 4& 4& 4&  2& 2& 2& 2& 2 \\
$10_{94}$:  &  &  &  &  &  &  &  & 2& 2& 2& 4& 4& 4& 6&  6&  8&  8& {\bf 8}   & 8 & 6 & 6 & 6 & 4& 4& 4& 2& 2& 2& 2& 2&  &  &  &  &  & \\
$10_{96}$:  &  &  & 2& 2& 2& 2& 4& 4& 4& 4& 6& 6& 6& 8&  8& 10& 10& {\bf 12}  & 10& 10& 8 & 8 & 6& 6& 4& 4& 4& 4& 2& 2& 2& 2&  &  &  & \\
$10_{110}$: &  &  &  &  & 2& 2& 2& 2& 4& 4& 4& 4& 6& 6&  8&  8& 10& {\bf 10}  & 10& 8&  8 & 6 & 6& 6& 4& 4& 4& 4& 2& 2& 2& 2&  &  &  & \\
$10_{112}$: &  &  &  & 2& 2& 2& 2& 2& 4& 4& 4& 6& 6& 8&  8&  8& 10& {\bf 10}  & 10& 10& 8 & 8 & 6& 6& 6& 4& 4& 4& 2& 2& 2& 2& 2&  &  & \\
$10_{117}$: & 2& 2& 2& 2& 2& 4& 4& 4& 4& 6& 6& 8& 8& 8& 10& 10& 12& {\bf{12}} & 12& 10& 10& 10& 8& 8& 6& 6& 6& 4& 4& 4& 4& 2& 2& 2& 2& 2
\end{tabular}
\end{small}

\subsection{Classification of alternating $10$-crossing knots with $u=1$}

Corollary~\ref{cor:Alt10}, together with known results (c.f.~\cite{Kawauchi})
suffice to classify all alternating $10$-crossing knots with $u=1$. 
Indeed, for alternating $10$-crossing knots, the previously known cases
can be reproved using only the Murasugi bound, cyclicity of 
$H_1(\Sigma(K);\Z)$, and Theorem~\ref{thm:AlternatingCriterion}. 

Specfically, consider those alternating, $10$-crossing knots with
$|\sigma(K)|\leq 2$ and cyclic $H_1(\Sigma(K);\Z)$, but which were not
covered by Corollary~\ref{cor:Alt10}.  Again, various knots fail
various of the tests in Theorem~\ref{thm:AlternatingCriterion}.

The following knots admit no even matching:
\[
10_{3}, 10_{19}, 10_{20}, 10_{24}, 10_{29}, 10_{36}, 10_{40}, 10_{65}, 10_{69}, 10_{89}, 10_{97}, 10_{108}, 10_{122}.
\]

The following knots have positive, even matchings, but none of these matchings is symmetric:
\[
10_{4}, 10_{11}, 10_{12}, 10_{13}, 10_{15}, 10_{16}, 10_{22}, 10_{28}, 10_{34}, 10_{35}, 10_{37}, 10_{38}, 10_{41}, 10_{43}, 10_{45}, 10_{115}.
\]

\section{Non-alternating $10$-crossing knots and the proof
of Corollary~\ref{cor:Rest10}}
\label{sec:NonAlternating}

The knots described in Corollary~\ref{cor:Rest10} are not
alternating. However, we claim that their branched double-covers are
$L$-spaces, and hence we can adapt the principles used in the proof of
Theorem~\ref{thm:AlternatingCriterion} (Montesinos' trick, followed by
Theorem~\ref{thm:LSpaceSymmetry}) only using correction terms for the
$\Sigma(K)$, in place of lengths of vectors of the Goeritz matrix.
The key problem remains to verify that $\Sigma(K)$ are $L$-spaces as
claimed, and then calculating the correction terms for $\Sigma(K)$. 

Some of these knots are Montesinos knots, and their branched
double-covers are Seifert fibered spaces. Hence the Heegaard Floer
homology can be calculated using techniques from~\cite{SomePlumbs}, as
explained in Subsection~\ref{subsec:Montesinos}. The remaining
cases are handled in Subsection~\ref{subsec:Remaining}.

\subsection{Corollary~\ref{cor:Rest10}: the Montesinos cases}
\label{subsec:Montesinos}

The knots in the list 
\begin{equation}
\label{eq:MontesinosList}
10_{125}, 10_{126}, 10_{130}, 10_{135}, 10_{138}. 
\end{equation}
are {\em Montesinos knots}, knots whose branched double-covers
are Seifert fibered spaces; in fact, the branched double-covers 
are the spaces with Seifert invariants 
\begin{equation}
\label{eq:SeifertFibereds}
\left(-2,\frac{1}{2},\frac{1}{3},\frac{4}{5}\right), 
\left(-2,\frac{1}{2}, \frac{2}{3}, \frac{1}{5}\right), 
\left(-2,\frac{1}{2},\frac{2}{3}, \frac{3}{7}\right), 
\left(-2,\frac{1}{2},\frac{1}{3}, \frac{2}{7}\right), 
\left(-2,\frac{1}{2},\frac{2}{5},\frac{2}{5}\right)
\end{equation}
respectively. Here our conventions on Seifert invariants are as
follows: $(b,\frac{\beta_1}{\alpha_1}, ..., \frac{\beta_n}{\alpha_n})$
are the Seifert invariants
for the three-manifold obtained as surgery on a configuration
consisting of a central circle, with surgery coefficient $b$, and a
collection of circles linking the central circle, with surgery
coefficients $-\frac{\alpha_i}{\beta_i}$.

For Seifert a fibered rational homology three-sphere $Y$, an algorithm
is given in ~\cite{SomePlumbs} which can be used to determine if $Y$ is
an $L$-space. Specifically, we start with a negative-definite plumbing
diagram for $\pm Y$, let $V$ denote the lattice generated by 
the vertices, and let $Q$ denote the induced bilinear form.
As presented, $V$ has a preferred basis given by the vertices of the
plumbing graph.
We consider the equivalence relation on
characteristic vectors in $V^*$ generated by
\begin{eqnarray*}
\kappa \sim \kappa \pm q(v) &&
{\text{if $v\in V$ is a preferred basis vector with 
$\pm\langle \kappa,v\rangle = Q(v,v)$}}.
\end{eqnarray*}
It is clear that if $\kappa\sim \lambda$
then $Q^*(\kappa\otimes \kappa)=Q^*(\lambda\otimes \lambda)$.
Let $X$ denote the number of equivalence classes with the property
that each $\kappa$ representing a given element of $X$ satisfies
the bound 
$$|\langle \kappa, v\rangle| \leq Q(v,v)$$
for each preferred basis vector $v\in V$. Clearly, the number of elements in
$X$ is at least as large as the number of elements in $\Coker q$.
In fact, according to~\cite{SomePlumbs}, the number of elements of $X$
is the rank of $\Ker U\subset \HFp(-Y)$, while it is elementary
to see that the number of elements
in the cokernel of $q$ is identified with the number of elements
in $H^2(Y;\Z)$. Thus, if the number of elements in $X$ agrees with
the number of elements in $H^2(Y;\Z)$, then $Y$ is an $L$-space.
Moreover, in~\cite{SomePlumbs}, it is shown that under these circumstances,
the map
$$\MinVecLen_Q\colon \Coker(q) \longrightarrow \Q$$
agrees with the correction terms for $-Y$
under a suitable identification of $\Coker(q)$ with $\SpinC(Y)$.

A straightforward calculation shows that all of the Seifert
fibered spaces in List~\eqref{eq:SeifertFibereds} satisfy this criterion,
and hence are $L$-spaces, and hence Theorem~\ref{thm:LSpaceSymmetry}
can be used to deduce the existence of a symmetry for some matching
of the vector of correction terms with our usual vector $B$.

For instance, for $\Sigma(10_{125})$, 
a matrix representing $Q$ in the preferred basis is
given by
\[
\left(\begin{array}{rrrrrrr}
-2& 1& 1& 1& 0& 0& 0 \\
1& -2& 0& 0& 0& 0& 0 \\
1& 0& -3& 0& 0& 0& 0 \\
1& 0& 0& -2& 1& 0& 0 \\
0& 0& 0& 1& -2& 1& 0 \\
0& 0& 0& 0& 1& -2& 1 \\
0& 0& 0& 0& 0& 1& -2
\end{array}\right).
\]
Calculating the function
$\MinVecLen_Q$, and ordering the elements of $\SpinC(Y)$
in a manner compatible with the group structure,
we get the vector $A$:
\[A= \left(\frac{1}{2},\frac{35}{22},\frac{19}{22},\frac{7}{22},
   - \frac{1}{22}  ,- \frac{5}{22}  ,
   - \frac{5}{22}  ,- \frac{1}{22}  ,\frac{7}{22},
   \frac{19}{22},\frac{35}{22}\right).
\]
Comparing against the corresponding vector $B$, we find that there are
no even, positive, symmetric matchings. Indeed, proceeding in a like manner
for all knots in List~\eqref{eq:MontesinosList},
we have that all even, positive matchings are given by
the following table:

\begin{tabular}{cccccccccccc}
$10_{125}$:&    &       &       &       &       &       {\bf2}& 2\\
$10_{126}$:&    &       &       &       &       &       {\bf2}& 2\\
$10_{130}$:&    &       &       &       2&      2&      {\bf2}& 2\\
$10_{135}$:&    &       2&      2&      2&      4&      {\bf4}& 4&      2&      2&      \\
$10_{135}'$:&   2&      2&      2&      2&      4&      {\bf4}& 4&      2&      2&      2\\
$10_{138}$:&    &       2&      2&      2&      4&      {\bf4}&         4&      2&      2&      2&      2
\end{tabular}

None of the above matchings is symmetric, and hence applying 
Montesinos' trick together with
Theorem~\ref{thm:LSpaceSymmetry} in the usual manner, we see that 
none of these knots has unknotting number equal to one.

\subsection{Corollary~\ref{cor:Rest10}: the remaining cases}
\label{subsec:Remaining}

The knots  $10_{148}$, $10_{151}$, $10_{158}$, and $10_{162}$ 
are not Montesinos knots. However,
their branched double-covers are $L$-spaces, and indeed we use
a refinement of Proposition~\ref{prop:Sharp}
to calculate the correction terms of $\Sigma(K)$.
Thus, we must construct sharp four-manifolds which bound the
knots in this list.

A knot projection of $K$ specifies a four-manifold 
$Z_K$
which bounds $\Sigma(K)$.
Starting from the white graph, we draw an unknot with
framing $0$ for all but one of the vertices in the white graph.  Next, we
associate to each edge in the white graph an unknot which links the
two $0$-framed unknots corresponding to the two vertices. This unknot is
given framing $+1$ if the crossing is consistent with the coloring
convention illustrated in Figure~\ref{fig:ColorConventions}, and it is
given framing $-1$ if it does not. This framing $\pm 1$ will be called
the sign of the edge. From $Z_K$ we can obtain another four-manifold $X_K$
by blowing down all of the unknots with
framing $\pm 1$ corresponding to the white graph. The remaining link
is obtained as a plumbing of unknots, with intersection form specified
by the Goeritz form for the projection of 
$K$. Specifically, the two-dimensional homology of $X_K$ is generated by
the vertices of the white graph, modulo the relation
$$\sum_{v\in V}v = 0,$$ 
where here $V$ denotes the vertices of the white graph.
Also, $Q(v\otimes v)$ is given by minus
the sum of the signs of all the edges leaving $v$, while if $v\neq w$,
then $Q(v\otimes w)$ is  the sum of the signs of all the edges
connecting $v$ and $w$. (In the alternating case, this construction
can be seen to be equivalent to the one given in
Section~\ref{sec:FirstApp}.)

\begin{figure}
\mbox{\vbox{\epsfbox{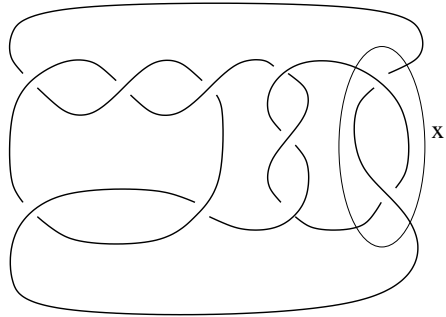}}}
\caption{\label{fig:10s148}
{\bf The knot $10_{148}$.}}
\end{figure}

\begin{figure}
\mbox{\vbox{\epsfbox{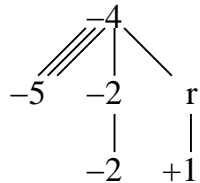}}}
\caption{\label{fig:Kirby10s148}
{\bf Kirby calculus pictured for $\Sigma(10_{148})$.}
	The integers here are framings on unknots, and the 
	lines represent linkings between the unknots. Setting $r=-1$,
	we obtain a picture of a four-manifold
	which bounds $\Sigma(10_{148})$.}
\end{figure}

For example, for the knot $K=10_{148}$ pictured in
Figure~\ref{fig:10s148}, we obtain a Kirby picture of $X_K$ which is a
plumbing of unknots, as pictured in Figure~\ref{fig:Kirby10s148},
after we blow down the circle labeled with $r=-1$. Ignoring this
circle (i.e. performing $r=\infty$ surgery), we obtain a picture for
the branched double-cover of $\Sigma(K_1)$, while setting $r=0$,
we get a picture for a branched double-cover of $\Sigma(K_0)$, 
where here $K_0$ and $K_1$ are the knots obtained by resolving
either of the intersection points in the oval marked by $x$ in 
Figure~\ref{fig:10s148}. Indeed, blowing down the unknot with
framing $+1$, we obtain a four-manifold with intersection form 
given by the matrix 
\[G=\left(
\begin{array}{rrrrr}
-4&3&1&0&1 \\
3&-5&0&0&0 \\
1&0&-2&1&0 \\
0&0&1&-2&0 \\
1&0&0&0& r-1
\end{array}
\right).
\]

In fact, let $X_2$ denote the four-manifold obtained by ignoring the
unknot with framing $r-1$ (whose intersection form is the top $4\times
4$ block of $G$). Now, $-\Sigma(K_0)$, $-\Sigma(K)$, and $-\Sigma(K_1)$
are related by a surgery long exact sequence, with cobordisms
\begin{eqnarray*}
W_0\colon -\Sigma(K_0)\longrightarrow -\Sigma(K)
&{\text{and}}&
W_1\colon -\Sigma(K) \longrightarrow -\Sigma(K_1).
\end{eqnarray*}
Moreover, if we blow down the obvious sphere of square $-1$ in the
composite cobordism $X_0=W_0\cup W_1 \cup X_2$, we obtain a four-manifold
$X_0'$ with intersection form given by $G$ with $r=0$, while
$X_1=W_1\cup X_2$ has intersection form given by $G$ with $r=-1$.  Our
aim is to show that $X_1$ is sharp.

To this end, we claim that the branched double-cover of $K_0$ is an
$L$-space, and indeed that the four-manifold with intersection form
given by $G$ with $r=-1$ is sharp. Indeed, the branched double-cover
of $\Sigma(K_0)$ is gotten by $+8$ surgery on the right-handed
trefoil, an $L$-space whose correction terms can be calculated for
example, using~\cite{SomePlumbs}; or
alternatively~\cite{AbsGraded}. Comparison with $M_Q$ for the
quadratic form for $G$ above with $r=0$, we conclude that the
four-manifold $X_0$ is sharp.

Moreover, $K_1$ is an alternating link and the top $4\times 4$
submatrix $G$ is the Goeritz matrix for an alternating
projection. Thus, $X_2$ is sharp.

Having verified that $-\Sigma(K_0)$ (a space with $8$ $\SpinC$
structures) and $-\Sigma(K_1)$ (a space with $23$ $\SpinC$ structures)
are $L$-spaces, it follows from Proposition~\ref{prop:LSpaces} that
$-\Sigma(K)$ (a space with $31$ $\SpinC$ structures) is an $L$-space
as well. Indeed, by Proposition~\ref{prop:Sharp}, it also now follows
that the matrix $G$ with $r=-1$ can be used to calculate the
correction terms $A$ for $\Sigma(K)$.

Again, there is a unique
even, positive matching, and it is given by
$$
\begin{array}{lc}
10_{148}: & 2, 2, {\bf 4}, 2, 2, 2.
\end{array}
$$

\begin{figure}
\mbox{\vbox{\epsfbox{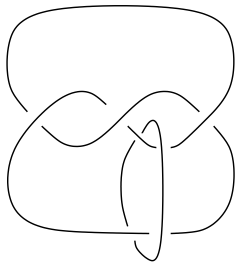}}}
\caption{\label{fig:QAltLink}
{\bf A link with determinant $8$.}
This is the link obtained by vertically resolving one of the crossings
in the oval $x$ for $10_{148}$ in Figure~\ref{fig:10s148}.}
\end{figure}

\begin{figure}
\mbox{\vbox{\epsfbox{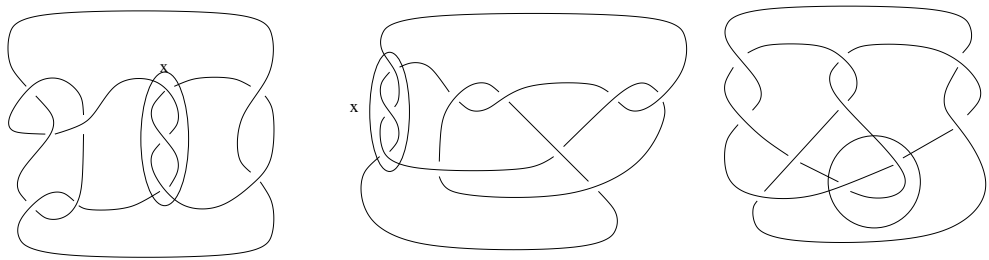}}}
\caption{\label{fig:NonAlt10}
{\bf The knots $10_{151}$, $10_{158}$, and $10_{162}$ respectively.}}
\end{figure}

We proceed similarly for the other three knots, which are pictured in
Figure~\ref{fig:NonAlt10}.

\begin{figure}
\mbox{\vbox{\epsfbox{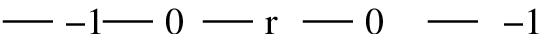}}}
\caption{\label{fig:Chain}
{\bf A chain of unknots.}
}
\end{figure}

For the knot $K=10_{151}$, the Kirby calculus description of $Z_K$
contains a certain chain unknots as illustrated in
Figure~\ref{fig:Chain} with $r=-1$, corresponding to the circled
region in the diagram for $10_{151}$. Blowing down the two $-1$-framed
and all the $+1$-framed ones unknots gives a four-manifold with
indefinite intersection form. Instead, we blow down all the
$+1$-framed unknots, and then perform
two handleslides and two handle cancellations,
after trading the $0$-framed unknots for one-handles. This replaces
the chain with a single unknot with framing $r-2$.

In this manner, we obtain a four-manifold $X(r)$ whose intersection form
is given by:
\[G=\left(\begin{array}{rrrrr}
-3&1&1&1&1 \\
1&-2&1&0&0 \\
1&1&-4&0&0 \\
1&0&0&r-2&0 \\
1&0&0&0&-2 \\
\end{array}
\right)\]
which bounds $\Sigma(K)$ when $r=-1$.

Setting $r=+1$, we obtain a four-manifold which bounds the branched
double-cover of the $(2,5)$ torus knot. This can be seen by replacing
the $-3/2$ twist insider the oval marked by $x$ in the picture for
$10_{151}$ in Figure~\ref{fig:NonAlt10} by a single $-1/2$ twist. A
direct calculation shows that this four-manifold is sharp. Moreover,
the $4\times 4$ submatrix obtained by deleting the last row and column
is a Goeritz matrix for the alternating knot with determinant $19$
obtained by forming the $K_1$-resolution of any of the three crossings
in the oval marked by $x$. In particular, the associated four-manifold
there is sharp, as well. Since the determinant of $G$ (for arbitrary
$r$) is given by $-24 + 19 r$, it follows from
Proposition~\ref{prop:Sharp} (and a descending induction
starting at $r=+1$)
that the four-manifold described by our plumbing of unknots with
framing $r\leq 1$ is sharp.

Thus, we can calculate the correction terms for $\Sigma(10_{151})$
using $G$ at $r=-1$ to get $A$. One then checks that the
only even, positive matching is given by
\begin{eqnarray*}
10_{151}: &&  2, 2, 2, 4, 4, {\bf 6}, 4, 4, 2, 2, 2, 2.
\end{eqnarray*}

We proceeding similarly for $K=10_{158}$. Again, we find a chain
is the plumbing description for $10_{158}$ as given in Figure~\ref{fig:Chain},
which we replace by an unknot with framing $r-2$.
Thus, we obtain a four-manifold $X_K$ which bounds
$\Sigma(K)$, with intersection form given by 
$$G=
\left(\begin{array}{rrrr}
-4&2&2&1 \\
2&-4&1&0 \\
2&1&-4&0 \\
1&0&0&r-2 \\
\end{array}
\right)$$
with $r=-1$. Ignoring the unknot with framing $r-2$,
we obtain the Goeritz matrix for the alternating link 
obtained from the $K_1$-resolution.
When $r=+1$, we obtain a still negative-definite
intersection form for a four-manifold which bounds
$\Sigma(4_1)$. Again, a direct calculation shows that this
four-manifold is sharp, and hence, by Proposition~\ref{prop:Sharp}, so
is $X_K$ (indeed since $G$ has matrix $25-20r$, the
corresponding four-manifold is sharp for all $r\leq 1$).  
Calculating the vector $A$ using the matrix $G$, we see
that there is a unique even, positive matching:
\begin{eqnarray*}
10_{158}: &&
2, 2, 2, 2, 4, 4, {\bf 6}, 4, 4, 2, 2, 2.
\end{eqnarray*}

Next, we turn to the knot $K=10_{162}$ (according to a numbering
scheme in which $10_{161}\neq 10_{162}$). Like in $10_{148}$, 
when constructing $X_K$, we can blow down a $+1$-framed unknot to
obtain a four-manifold $X_K$ with intersection form given by
$$
G=\left(
\begin{array}{rrrr}
-5&2&2&0 \\
2&-5&2&0 \\
2&2&-4&1 \\
0&0&1&r-1 \\
\end{array}
\right)$$
when $r=-1$.  The $3\times 3$ submatrix on the upper left corner
specifies an intersection form for the manifold which bounds
$\Sigma(K_1)$, while setting $r=0$ specifies an intersection form for
a four-manifold which bounds $\Sigma(K_0)$.  It is easy to see that
both $K_1$ and $K_0$ have alternating projections -- $K_1$ is an
$8$-crossing link with determinant $28$, while $K_0$ is the knot $5_2$
(with determinant $7$). Comparing against the Goeritz matrix for the
alternating projection, we see that the four-manifold specified 
by $G$ with
$r=0$ is sharp (for $\Sigma(K_0)$.
Thus, by Proposition~\ref{prop:Sharp}, the four-manifold $X_K$
is sharp, as well and hence the correction terms for $\Sigma(10_{162})$
can be calculated using
the matrix $G$ with $r=-1$.

The unique even, positive matching is given now by
\begin{eqnarray*}
10_{162}: &&  2, 2, 4, {\bf 4}, 4, 2, 2, 2.
\end{eqnarray*}

Note that none of the even, symmetric matchings for $10_{148}$,
$10_{151}$, $10_{158}$, and $10_{162}$ listed above is symmetric, and
hence none of these knots has unknotting number equal to one, by
Theorem~\ref{thm:LSpaceSymmetry}.

\subsection{Final remarks}

\begin{figure}
\mbox{\vbox{\epsfbox{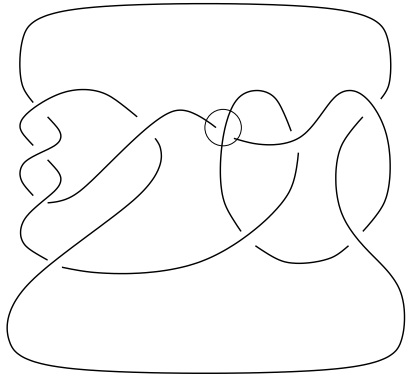}}}
\caption{\label{fig:10s131}
{\bf Unknotting $10_{131}$.}}
\end{figure}

In fact, the methods descibed here -- cyclicity of $H_1(\Sigma(K);\Z)$,
and the methods from Theorem~\ref{thm:AlternatingCriterion} --
are sufficient to classify all $10$-crossing knots with $u=1$, with 
the two exceptions $10_{145}$ and $10_{153}$. In particular, the knots
$10_{140}$, $10_{163}$, $10_{165}$ admit no even matching,
while $10_{144}$ admits no even, positive matching.
(Note that $10_{131}$ is listed as
having unknown unknotting number in some sources, but it can be
unknotted in one step according to Figure~\ref{fig:10s131}, by
changing the indicated crossing, see also~\cite{Stoimenow}.)

The knot $10_{145}$ has $u=2$ according to~\cite{Tanaka}, while
$10_{153}$ has $u=2$ according to~\cite{GordonLueckeUK1}, completing
the classification of all knots with $u=1$ and $\leq 10$ crossings.

%       \subsection{Final things}
%        
%       Discuss (previously known)
%       \[
%        10_{140},  10_{144},  10_{162},  10_{165}.
%       \]

\section{Refinements}
\label{sec:Refine}

\subsection{Signed unknottings}

The reason for the two choices of $\epsilon$ in the statements of
Theorems~\ref{thm:AlternatingCriterion} rests on orientations. On the
one hand, the condition on a knot of having unknotting number one does
not distinguish a knot $K$ from its mirror $r(K)$. By contrast, the
branched double-cover of $S^3$ with its standard orientation along a
knot $K$ endows $\Sigma(K)$ with a natural orientation, with the
property that $\Sigma(r(K))\cong -\Sigma(K)$; moreover the signs
of the correction terms depend on a choice of orientation for its underlying
three-manifold. 

\begin{figure}
\mbox{\vbox{\epsfbox{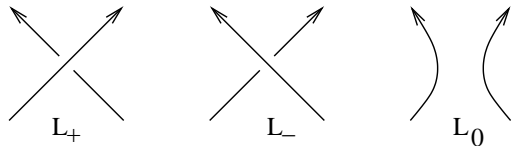}}}
\caption{\label{fig:SignConventions}
{\bf Sign conventions on crossings.} Crossings of type $L_+$ are
positive, $L_-$ is negative, and $L_0$ is the oriented resolution.}
\end{figure}

A refined statement of this result can be formulated which makes use
of a more orientation-dependent notion of unknotting number one: we
could consider knots $K$  which can be unknotted by changing a
single negative crossing to a positive crossing for some projection
$K$ (we use here the usual conventions from knot theory as illustrated in
Figure~\ref{fig:SignConventions}). One obstruction to this sign-refined
question is the signature $\sigma(K)$ of the knot. Specifically, recall that
the signature
satisfies the  inequality
\begin{equation}
\label{eq:SigmaIneq}
\sigma(K_-)-2\leq \sigma(K_+)\leq \sigma(K_-),
\end{equation}
and hence a knot with this property has $\sigma=0$ or $2$.

For the purpose of the following statement, recall that we use
the coloring convention illustrated in 
Figure~\ref{fig:ColorConventions}.
Given an alternating knot with determinant $D$, 
fix a regular alternating projection
and fix a corresponding Goeritz matrix $Q$, using the white graph
as described in the introduction.
Given an isomorphism $\phi\colon\Zmod{D}\longrightarrow \Coker(q)$, 
let
\[
T_{\phi}(i)=(-1)^{\frac{\sigma(K)}{2}}\cm \MinVecLen_Q(\phi(i))-\gamma_{D}(i). 
\]

\begin{theorem}
\label{thm:SignRefinedVersion}
Let $K$ be an alternating knot which can be unknotted by changing a
single negative crossing to a negative one in some (not necessarily
alternating) projection of $K$,
then there is an isomorphism $\phi\colon
\Zmod{D}\longrightarrow \Coker(q)$
with the properties that for all $i\in\Zmod{D}$:
\begin{eqnarray*}
T_{\phi}(i)&\equiv& 0 \pmod{2}
\\
T_{\phi}(i)&\geq& 0.
\end{eqnarray*}
If in addition $|\MinVecLen_{Q_K}(0)|\leq \OneHalf$, then
there is a choice of $\phi$ satisfying the above to constraints, and
the following additional symmetry:

\[
T_{\phi}(i)=T_{\phi}(2k-i)
\]
for $1\leq i < k$ when $D=4k-1$ and 
for $0\leq i < k$ when $D=4k+1$.
\end{theorem}

\begin{proof}
We need the following precise version of Montesinos' lemma,
see~\cite{Stoimenow}: if $K$ is a knot with determinant $D=2n+1$ which
can be unknotted by changing a negative crossing to a positive one,
then $\Sigma(K)=S^3_{-\epsilon\cm\frac{2n+1}{2}}(C)$ where
$\epsilon=(-1)^{\frac{\sigma}{2}}$. (Note that if $K$ can be unknotted
by changing a negative crossing to a positive one, then $\sigma(K)=0$
or $2$ according to Inequality~\eqref{eq:SigmaIneq}.)

We have the two triples of three-manifolds $(\Sigma(K_0), \Sigma(K_+), \Sigma(K_1)$ and $(\Sigma(K_1), \Sigma(K_-), \Sigma(K_0))$
which are related by two-handle additions as in the hypothesis of
Theorem~\ref{thm:ExactSeq}. Let $A$ denote the two-handle 
cobordism from $\Sigma(K_1)$ to $\Sigma(K_-)$ and $B$ denote the 
two-handle from $\Sigma(K_+)$ to $\Sigma(K_1)$.
Note that here $K_0$ is the oriented resolution of $K_{\pm}$.
By handlesliding, it is easy to see that $A\circ B$ contains a sphere
with square $-2$, and another linking two-handle. 
Our assumption that $K_+$ is the unknot then ensures that $\Sigma(K_-)$
can be written as $S^3_{\pm \frac{2n+1}{2}}$, where the sign
depends on $b_2^+(A\circ B)$. 

It is a standard fact that for a knot $K$
$(-1)^{\frac{\sigma(K)}{2}}D=\Delta_{K}(-1)$ (c.f.~\cite{Lickorish}).
Suppose that $\sigma(K)$ is zero. Then the
above fact, together with 
the skein relation for the Alexander polynomial gives that
$\det(K_0)=n$. It follows that $\det(K_1)=n+1$. These
in turn ensure that both cobordisms $A$ and $B$ are negative-definite,
and hence that $\Sigma(K)=\Sigma(K_-)=S^3_{-\frac{2n+1}{2}}(C)$.

In the case where $\sigma(K)=2$, the same argument now proves that
$\det(K_0)=n+1$ and $\det(K_1)=n$, from which it follows that
$b_2^+(A)=1$. Thus, the Kirby calculus picture for $A\circ B$ consists
of a knot $C$ with positive framing, with a linking unknot with
framing $-2$. By modifying the cobordism $A\circ B$ in a
straightforward way, we can trade the linking unknot for another one
with framing $+2$, at the cost of increasing the framing on $C$ by
one. Thus, we have expressed $\Sigma(K)=\Sigma(K_-)\cong
S^3_{\frac{2n+1}{2}}(C)$.

With the signs pinned down, now, the proof of the result follows from
the proof of Theorem~\ref{thm:AlternatingCriterion}.
\end{proof}

More informally, if an alternating knot has unknotting number equal to
one and $|\MinVecLen_{Q_K}(0)|\leq \OneHalf$, but it has positive,
even, symmetric matchings with only one choice of $\epsilon$, then 
it does not have two one-step unknottings with different signs.

Some knots -- for example, $8_{13}$ -- can be unknotted by single
crossing-changes with either sign. According to
Inequality~\eqref{eq:SigmaIneq}, a knot with this property must have
vanishing signature.  Thus, we illustrate
Theorem~\ref{thm:SignRefinedVersion} with a knot whose signature
vanishes.

\begin{figure}
\mbox{\vbox{\epsfbox{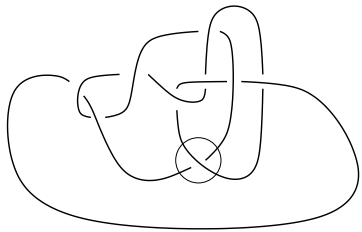}}}
\caption{\label{fig:9s33}
{\bf The knot $9_{33}$.}
	By changing the circled indicated (negative) crossing,
	we can unknot $9_{33}$. However, no projection of 
	$K=9_{33}$ has a positive crossing which, when changed,
	unknots $K$.}
\end{figure}

Consider the knot $9_{33}$ illustrated in Figure~\ref{fig:9s33}. For
one choice of $\epsilon$, there is a single even, positive matching of
the form $$ 2, 2, 2, 2, 4, 4, 4, 6, 6, {\bf 8}, 6, 6, 4, 4, 2, 2, 2,
2,$$ which is evidently asymmetric, while for the opposite choice of
$\epsilon$, we have the following unique even, positive matching which
is symmetric:
\begin{equation}
\label{eq:AllowedMatching}
2, 2, 2, 2, 4, 4, 4, 6, 6, {\bf 8}, 6, 6, 4, 4, 4, 2, 2, 2, 2
\end{equation}
(the existence of
at least one such a matching is guaranteed
from Theorem~\ref{thm:AlternatingCriterion}, as 
$9_{33}$ has unknotting number equal to one). 
This shows, however, that there is no one-step unknotting which
involves a crossing change with the opposite sign.

\subsection{Interpreting the matchings}

If $K$ is an alternating knot with unknotting number equal to one
and $|\MinVecLen_{Q_K}(0)|\leq \OneHalf$, then the
even integers $T_{\phi,\epsilon}$ guaranteed by Theorem~\ref{thm:AlternatingCriterion} 
have a concrete topological interpretation, which we now explain.

If $K$ is a knot with unknotting number one, then of course we can
draw an arc $\gamma$ in $S^3$ which connects two points of $K$, with
the property that a standard modification of $K$ in a tubular
neighborhood of $\gamma$ gives us the unknot. From a dual point of
view, a knot with unknotting number equal to one can be 
specified by an unknot together with a (framed) arc $\delta$ (in a
neighborhood of which we modify the unknot to get $K$). We call $\delta$
a {\em knotting arc} for $K$.

On the other hand, if $\delta$ is an arc connecting two points on an
unknot $O$, we can construct a knot $C$ in $S^3$, thought of as the
branched double-cover of $O$.  We claim that if $K$ is a knot which
satisfies all the hypotheses of
Theorem~\ref{thm:AlternatingCriterion}, then the (even) integers
appearing as differences between correction specify the Alexander
polynomial of $C$ (up to some finite indeterminacy determined by the
possible choices of $\epsilon$ and $\phi$).

To state the result, it is convenient to reformulate the information in the Alexander polynomial of $C$.
Let $C$ be a knot in $S^3$, and write its symmetrized Alexander polynomial as
$$\Delta_C(T)=a_0+\sum_{i>0} a_i(T^i+T^{-i}).$$
then, its  torsion coefficients $\MT_i(C)$ are given by the formula:
$$\MT_i(C)=\sum_{j=1}^\infty j\cm a_{|i|+j}.$$

It will also be convenient to have the following notation.
If $C\subset S^3$ is a knot, then an orientation for $C$ specifies a map
$$\sigma\colon \Zmod{p}\longrightarrow S^3_p(C)$$
by the condition that $\sigma(i)$ extends over the two-handle cobordism
$W_p(C)$ from $S^3$ to $S^3_p(C)$ as a $\SpinC$ structure $\spinc$ with
$$\langle c_1(\spinc),[{\widehat C}] \rangle \equiv 2i-p\pmod{2p}.$$

\begin{theorem}
\label{thm:LSpaceAlex}
Let $C\subset S^3$ be a knot in $S^3$ with the property that $S^3_p(C)$ is an $L$-space $Y$, then
$$2\MT_i(C)=
\left\{\begin{array}{ll}
-d(S^3_p(C),\sigma(i)) + d(S^3_p(O),\sigma(i)) &{\text{if $2|i|\leq p$}} \\
0 & {\text{otherwise,}}
\end{array}
\right. $$
\end{theorem}

The above is essentially a restatement of
Corollary~\ref{AbsGraded:cor:AlexLens} of~\cite{AbsGraded}, only that
result is stated in the case where $Y$ is the lens space $L(p,q)$;
however, the only property about lens spaces used in its proof is that
lens spaces are $L$-spaces. (In fact, the result is seen as a consequence
of a stronger result, which describes $\HFp(S^3_0(C))$ in terms of 
the correction terms for $S^3_p(C)$ and the map $\sigma$.)

Thus, in view of Theorem~\ref{thm:LSpaceFracSurg}, if $K$ is an
alternating knot with unknotting number equal to one, then twice
the torsion coefficients of the branched double-cover of the knotting
arc must appear (in order) in some matching for $K$. Incidentally, 
from this point of view, the symmetric condition on the matching
corresponds to the usual symmetry of the Alexander polynomial of $C$.

For example, we saw earlier that the knot $9_{33}$ has a unique even,
positive, and symmetric matching, as given in
Equation~\eqref{eq:AllowedMatching}. Indeed, converting from the
torsion back to the Alexander polynomial, it follows that the
Alexander polynomial of $C$ is given by
$$
-1+(T^2+T^{-2})-(T^4+T^{-4})+(T^5+T^{-5})-(T^8+T^{-8})+(T^9+T^{-9})
$$
which, incidentally, is the Alexander polynomial of the $(7,4)$ torus knot.
It is, in fact, reasonable to expect that if $\delta$ is any knotting arc
of $9_{33}$, then its branched double-cover is the $(7,4)$ torus knot.
(Compare Berge's conjecture on knots which admit
lens space surgeries, c.f.~\cite{Berge}, \cite{Kirby}.)

\subsection{Alternating knots with unknotting number one as a source of examples}

Alternating knots with unknotting number equal to one can be viewed as
a wide source of knots in $S^3$ which admit $L$-space surgeries
(by taking the branched double-cover of the knotting arc).  Such
knots are very special. For example, in~\cite{NoteLens}, we prove the following:

\begin{theorem}
\label{thm:StructAlex}
Suppose that $C$ is a knot in $S^3$ which admits an integral $L$-space
surgery, then all the coefficients of its Alexander polynomial are
$\pm 1$ and the non-zero coefficients all alternate in sign.
\end{theorem}

The above theorem appears as Corollary~\ref{NoteLens:cor:StructAlex}
in~\cite{NoteLens}, where it is seen as a corollary to a more general
result which constrains the structure of the ``knot Floer homology'' of
$K$, c.f.~\cite{HolDiskKnots} and~\cite{RasmussenThesis}.
Indeed combining results from~\cite{NoteLens},~\cite{GenusBounds}
and~\cite{4BallGenus}, we get that for such a knot, the Seifert genus,
the four-ball genus, and the degree of the Alexander polynomial all
agree.  Moreover, hyperbolic knots with $L$-space surgeries all
provide infinitely many examples of hyperbolic three-manifolds which
admit no taut foliation~\cite{GenusBounds}.  (Compare also~\cite{KMOS}
for analogous results in the realm of Seiberg-Witten monopole Floer
homology.)

Berge's construction~\cite{Berge} gives many examples knots with
$L$-space surgeries. Alternating knots with unknotting number equal to
one provide another source of such examples.

\subsection{Stronger forms of Theorem~\ref{thm:AlternatingCriterion} using
  $L$-space surgeries theorems}

Results from~\cite{NoteLens} stated above concerning the structure of
the Alexander polynomial of a knot admitting $L$-space surgeries can
be viewed as giving further restrictions on the positive, even,
symmetric matchings associated to alternating knots with unknotting
number equal to one.  For example, we obtain the following result:

\begin{theorem}
	\label{thm:StrongerForm}
  If $K$ is an alternating knot with unknotting number equal to one
  with determinant $D=4k\pm 1$, then there is a choice of isomorphism
  $\phi\colon \Zmod{D}\longrightarrow \Coker{q}$ and $\epsilon$, with
  the property that the matching $T_{\phi,\epsilon}$ satisfies the following
  restrictions:
  \begin{itemize}
     \item $T_{\phi,\epsilon}(i)\equiv 0 \pmod{2},$
     \item $T_{\phi,\epsilon}(i)\geq 0$ for all $i$;
       moreover for 
       for $i=1,...,k-1$
       \begin{equation}
	\label{eq:TorsionIncreases}
	T_{\phi,\epsilon}(i)\leq T_{\phi,\epsilon}(i+1)
       \leq T_{\phi,\epsilon}(i)+2
\	\end{equation}
     \item $T_{\phi,\epsilon}(i)=T_{\phi,\epsilon}(2k-i)$ for $1\leq
     i<k$.  If, in addition, $|\MinVecLen_{Q_K}(0)|\leq \OneHalf$, the
     above symmetry extends to $i=0$ when $D=4k+1$.  \end{itemize}
\end{theorem}

The proof is based on a combination of Theorem~\ref{thm:StructAlex},
together with the techniques of this paper. Before giving the proof,
we give the following improvement of Theorem~\ref{thm:LSpaceFracSurg}
(which we state in the notation from
Section~\ref{sec:LSpaceSurgeries}):

\begin{theorem}
\label{thm:LSpaceFracSurgStrong}
Assume that $C\subset S^3$ is a knot with the property that for some $n>1$
$S^3_{-(2n-1)/2}(C)$ is an $L$-space. Then so is $S^3_{-n}(C)$, and 
also for
$i=1,...,2n-2$,
$$d(S^3_{-n}(C),v_i) - d(S^3_{-n}(O),v_i)
= d(S^3_{-(2n-1)/2}(C),w_i)-d(S^3_{-(2n-1)/2}(O),w_i).$$
\end{theorem}

\begin{proof}
  We proceed as in the proof of Theorem~\ref{thm:LSpaceFracSurg}, only
  we no longer have a hypothesis which ensures that
  $d(S^3_{-n}(C),v_j)-d(S^3_{-n}(O),v_j)=0$ for $v_j=v_0$ or $v_{2k}$
  (depending on the parity of $n$ as in that proof); we have only that
  $0\leq d(S^3_{-n}(C),v_j)-d(S^3_{-n}(0),v_j)$. However, by
  Theorem~\ref{thm:StructAlex}, it follows that
  $d(S^3_{-n}(C),v_j)-d(S^3_{-n}(O),v_j)\leq 2$ (since this difference
  is twice some coefficient of the Alexander polynomial of $C$). It
  follows now that from this, together with the argument from the
  proof of Theorem~\ref{thm:LSpaceFracSurg} that the inequality
  $$d(S^3_{-n}(C),v_i)-d(S^3_{-n}(O),v_i) \leq
  d(S^3_{-\frac{2n-1}{2}}(C),w_i)-d(S^3_{-\frac{2n-1}{2}}(O),w_i)$$
  (from Lemma~\ref{lemma:Positivity})
  can be a strict inequality at most one value of $i=0,...,2n-2$ (note
  that $d(S^3_{-n}(C),v_i)-d(S^3_{-n}(O),v_i)\equiv 0\pmod{2}$).  By
  the symmetry of the $w_i$ sending $i\mapsto 2n-1-i$, it follows at
  once that the inequality can be strict only when $i=0$.
\end{proof}

\vskip.3cm
{\noindent{\bf Proof of Theorem~\ref{thm:StrongerForm}.}}  Proceeding
as in the earlier proof, and applying
Theorem~\ref{thm:LSpaceFracSurgStrong} if necessary, we see that for
$0\leq i<k$ the numbers
$T_{\phi,\epsilon}(k-i)$ are differences in correction terms
for $S^3_{-n}(C)$ and $S^3_{-n}(O)$. In turn, according to 
Theorem~\ref{thm:LSpaceAlex}, these are identified with torsion coefficients
for $C\subset S^3$. Now, Equation~\eqref{eq:TorsionIncreases} is
a consequence of this fact, together with Theorem~\ref{thm:StructAlex}.
\qed

As we saw, Theorem~\ref{thm:AlternatingCriterion} suffices for the
study of knots with $\leq 10$, but it is possible that for other
applications, the stronger form given in
Theorem~\ref{thm:StrongerForm} might be useful.

\commentable{
\bibliographystyle{plain}
\bibliography{biblio}
}

\end{document}